\documentclass[fontsize=11pt,a4paper,DIV11]{scrartcl}
\usepackage[utf8]{inputenc} 
\usepackage{url}

\usepackage{mathtools}

\usepackage{graphicx}
\usepackage{refcount}
\def\NI{\noindent}

\def\sk{\smallskip}

\def\bs{\bigskip}

\def\Proof{\par{\NI\bfseries Proof.~}}
\def\ProofOf#1{\par{\NI\bfseries Proof of #1.}~}
\def\qed{\hfill\fbox{\hbox{}}\bs}

\usepackage{amssymb}
\newtheorem{theorem}{Theorem}
\newtheorem{lemma}{Lemma}

\newtheorem{proposition}{Proposition}

\newtheorem{corollary}{Corollary}

\def\Case#1.{\rih{Case #1.}}
\def\Phase#1.{\NI{\bfseries Phase #1.}}
\def\Fact#1.{\par\sk{\NI\bfseries Fact #1.}}

\def\ITEMMACRO #1 ??? #2 ???{\par\medskip\noindent%
\hangindent=#2em\setbox0\hbox{#1 \kern5pt}%
\ifdim\wd0<\hangindent\setbox0\hbox to\hangindent{\hss#1\quad}\fi%
\box0\ignorespaces}

\graphicspath{{Figures/}}
\input{psfig.sty}
\usepackage{subfig}
\usepackage{color}
\usepackage{tikz}
\usetikzlibrary{arrows,shapes}
\usetikzlibrary{decorations.markings}
\tikzstyle{edge} = [draw,thick,-]
\tikzstyle directed=[postaction={decorate,decoration={markings,
			    mark=at position .65 with {\arrow[arrowstyle]{stealth}}}}]

\DeclareMathOperator{\idim}{idim}
\newcommand{\dom}{\text{prod}}
\newcommand{\RR}{\mathbb{R}}

\newcommand{\Min}{\textup{Min}}
\newcommand{\Max}{\textup{Max}}
\usepackage{bibitemsep}

\title{Grid Intersection Graphs and Order Dimension\thanks{Steven Chaplick is supported by ESF EuroGIGA project GraDR,
        Stefan Felsner is partially supported by DFG grant FE-340/7-2 and ESF EuroGIGA project GraDR
        , Udo Hoffmann and Veit Wiechert are supported by the Deutsche Forschungsgemeinschaft within the
        research training group 'Methods for Discrete Structures' (GRK 1408) 
  }}

\author{
{\scshape Steven Chaplick}
{and} {\scshape Stefan Felsner}
{and}\\[0pt]
 {\scshape Udo~Hoffmann}
{and} {\scshape Veit Wiechert}\\[1pt]
\small
\url{{chaplick,felsner,uhoffman,wiechert}@math.tu-berlin.de}\\
\small
Institut f\"ur Mathematik,
         Technische Universit\"at Berlin\\[-4pt]
\small
Stra\ss e des 17. Juni 136,
         D-10623 Berlin, Germany
}

\date{24. November 2015}
\begin{document}
\maketitle

\begin{abstract}
\noindent
We study subclasses of grid intersection graphs from the perspective of order
dimension.  We show that partial orders of height two whose comparability
graph is a grid intersection graph have order dimension at most four. Starting
from this observation we provide a comprehensive study of classes of graphs
between grid intersection graphs and bipartite permutation graphs and the
containment relation on these classes. Order dimension plays a role in many
arguments.
\end{abstract}

% -------------------------------------------------------------
\section{Introduction}
% -------------------------------------------------------------
One of the most general standard classes of geometric intersection graphs is the
class of \emph{string graphs}, i.e., the intersection graphs of curves in the plane. 
String graphs were introduced to study electrical networks~\cite{sinden1966topology}.
The \emph{segment intersection graphs} form a natural subclass of string graphs, where
the curves are restricted to straight line segments.  We study subclasses
where the line segments are restricted to only two different slopes and
parallel line segments do not intersect. This class is known as the class
of \emph{grid intersection graphs} (GIG). An important feature of this class
is that the graphs are bipartite.
Subclasses of GIGs appear in several technical applications.  For example in
\emph{nano PLA-design}~\cite{shrestha2011two} and for detecting \emph{loss of
  heterozygosity events in the human
  genome}~\cite{halldorsson11}.

Other restrictions on the geometry of the representation are used to study
algorithmic problems.  For example, \emph{stabbability} has been used to study
hitting sets and independent sets in families of
rectangles~\cite{chepoi2013approximating}.
Additionally, computing the jump number of a poset, which is NP-hard in
general, has been shown solvable in polynomial time for bipartite posets with
interval dimension two using their restricted GIG
representation~\cite{soto2011jump}.

Beyond these graph classes that have been motivated by applications and
algorithmic considerations, we also study several other natural intermediate
graph classes.  All these graph classes and properties are formally defined in
Subsection~\ref{subsec:def}.

\begin{figure}[htb]
    \centerline{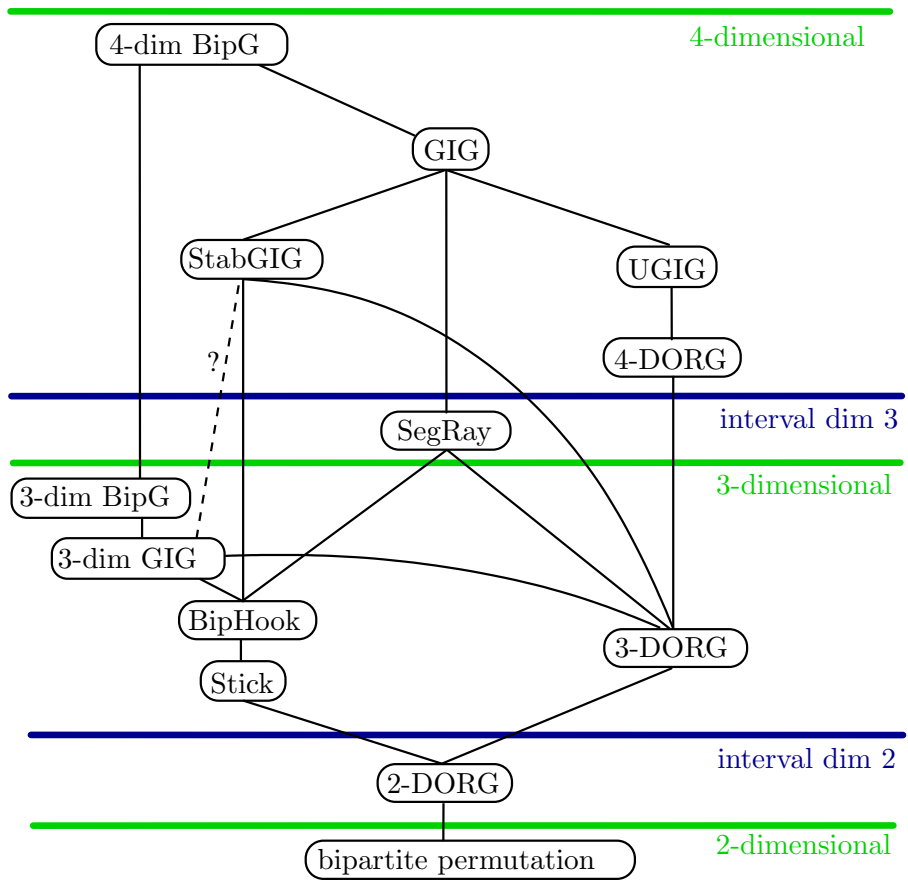}
    \caption{The inclusion order of graph classes studied in this paper.\label{fig:classinclusion}}
    \end{figure}

The main contribution of this work is to establish the strict containment and
incomparability relations depicted in Figure~\ref{fig:classinclusion}.
We additionally relate these classes to incidence posets of planar and outerplanar graphs.

In Section~\ref{sec:relations} we use the geometric representations to
establish the containment relations between the graph classes as shown in
Figure~\ref{fig:classinclusion}.  The maximal dimension of graphs
in these classes is the topic of Section~\ref{sec:dimension}.  In
Section~\ref{sec:incidence} we use vertex-edge incidence posets of planar graphs to
separate some of these classes from each other.  Specifically, we show that
the vertex-edge incidence posets of planar graphs are a subclass of \emph{stabbable} GIG
(StabGIG), and that vertex-edge incidence posets of outerplanar graphs are a subclass of
\emph{stick intersection graphs} (Stick) and \emph{unit} GIG (UGIG).  The
remaining classes are separated in Section~\ref{sec:separating}.
The separating examples are listed in Table~\ref{tab:separating}.
As part of this, we show that the vertex-face incidence posets of outerplanar graphs are
\emph{segment-ray intersection graphs} (SegRay). As a corollary we
obtain that they have interval dimension at most $3$.

%%%%%%%%%%%%%%%%%%%%%%%%%%%%%%%%%%%%%%%%%%%%%%%%%%%%%%%%%%%%%%%%%%%%%%%%%%
% in einem figure environment mit caption
\calc_figscale{88}%

%%%%%%%%%%%%%%%%%%%%%%%%%%%%%%%%%%%%%%%%%%%%%%%%%%%%%%%%%%%%%%%%%%%%%%%%%%

\begin{table}
\begin{center}\small
\begin{tabular}{p{4cm}p{3cm}p{33mm}}
\hline
  Class I  \hfill $\not\subseteq$    & Class II                & Example \\[2pt] \hline \hline
  GIG                          & 3-dim BipG               & $S_4$ \\ \hline
  3-dim BipG                    & GIG/3-dim GIG           & Proposition~\ref{prop:threedimnonpseudoseg} \\ \hline
  3-dim GIG                    & SegRay                  & $P_{K_4}$, Proposition~\ref{prop:segrayouterplanar}\\
                               & 3-dim BipG               & Proposition~\ref{prop:threedimnonpseudoseg} \\ \hline
  StabGIG                      & 3-dim GIG               & $S_4$\\
                               & SegRay                  & $S_4$\\ \hline
  SegRay                       & 3-dim GIG               & Proposition~\ref{prop:segrayoutermap}\\
                               & StabGIG                 & Proposition~\ref{prop:nonstabsegray}\\ \hline
  UGIG                         & 3-dim GIG               & $S_4$\\
                               & StabGIG                 & Proposition~\ref{prop:4dorgnotstabbable}\\
                               & 4-DORG                  & $C_{14}$, see~\cite{STU10}\\
                               & SegRay                  & $S_4$, Proposition~\ref{prop:segrayidim}\\ \hline
  BipHook                      & 3-DORG                  & Trees \\
                               & Stick                   & Proposition~\ref{prop:hooknotstick}\\ \hline
  Stick                        & UGIG                    & Proposition~\ref{prop:STICKnotUGIG}\\ \hline
  4-DORG                       & 3-dim GIG               & $S_4$\\
                               & StabGIG                 & Proposition~\ref{prop:4dorgnotstabbable}\\
                               & SegRay                  & $S_4$\\
                               & 3-DORG                  & $S_4$\\ \hline
  3-DORG                       & BipHook                 & Proposition~\ref{prop:3dorgnothook}\\ \hline
  2-DORG                       & 2-dim BipG               & $S_3$\\ \hline
\end{tabular}
\end{center}
\caption{Examples separating graph classes in
Figure~\ref{fig:classinclusion}\label{tab:separating}}
\end{table}

%-------------------------------------------------------------
\subsection{Definitions of Graph Classes}\label{subsec:def}
%-------------------------------------------------------------
%%%%%%%%%%%%%%%%%%%%
% in einem figure environment mit caption
\calc_figscale{30}%
\begin{figure}[htb]
    \centerline{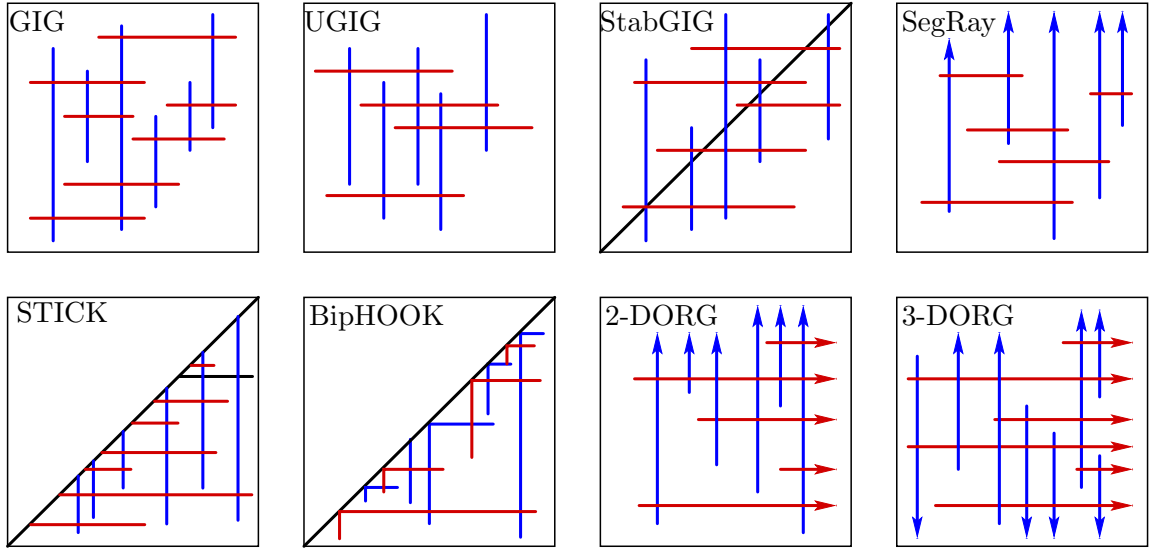}
    \caption{Typical intersection representations of graphs 
in the graph classes studied in this paper.\label{fig:examples}}
    \end{figure}%

%%%%%%%%%%%%%%%%%%%%%%%%%%%%%%%%%%%%%%%%%%%%%%%%%%%%%%%%%%%%%%%%%%%%%%%%%%

We introduce the graph classes from Figure~\ref{fig:classinclusion}.
A typical drawing of a representation is shown in Figure~\ref{fig:examples}.
We denote the class of bipartite graphs by BipG.
A \emph{grid intersection graph} (GIG) is an intersection graph of horizontal
and vertical segments in the plane where parallel segments do not intersect.
Some authors refer to this class as \emph{pure-GIG}.
If $G$ admits a grid intersection representation such that all segments have
the same length, then $G$ is a \emph{unit grid intersection graph} (UGIG).

A segment $s$ in the plane is \emph{stabbed} by a line $\ell$ if $s$ and
$\ell$ intersect.  A graph $G$ is a \emph{stabbable grid intersection graph}
(StabGIG) if it admits a grid intersection representation such that there
exists a line that stabs all the segments of the representation.  Stabbable
representations are generally useful in algorithmic settings as they provide a
linear ordering on the objects involved,
see~\cite{chepoi2013approximating,soto2014independent}.

A \emph{hook} is the union of a horizontal and a vertical segment that share
the left respectively top endpoint. The common endpoint, i.e., the bend of the hook,
is called the \emph{center} of the hook.  A graph $G$ is a \emph{hook graph}
if it is the intersection graph of a family of hooks whose centers are all on
a line $\ell$ with positive slope (usually $\ell$ is assumed to be the line
$x=y$). Hook graphs have been introduced and studied
in~\cite{maxpointtolerance,halldorsson11}, \cite{tomMaster},
and~\cite{SotoThraves2015}.  The graphs are called \emph{max point-tolerance graphs} in ~\cite{maxpointtolerance} and \emph{loss of heterozygosity graphs} in~\cite{halldorsson11}.  Typically these graphs are not
bipartite. We study the subclass of bipartite hook graphs
(BipHook).

A hook graph admitting a representation where every hook is degenerate, i.e.,
it is a line segment, is a \emph{stick intersection graph} (Stick).  In other words, Stick
graphs are the intersection graphs of horizontal or vertical segments that
have their left respectively top endpoint on a line $\ell$ with positive slope.

Intersection graphs of rays (or half-lines) in the plane have been previously
studied in the context of their chromatic number~\cite{kostochka1998coloring}
and the clique problem~\cite{cabello2013clique}.  We consider some natural
bipartite subclasses of this class.  Consider a set of axis-aligned rays in
the plane. If the rays are restricted into two orthogonal directions, e.g. up
and right, their intersection graph is called a \emph{two directional orthogonal
  ray graph} (2-DORG). This class has been studied in~\cite{STU10} and
\cite{soto2011jump}. Analogously, if three or four directions are allowed for
the rays, we talk about 3-DORGs or 4-DORGs. The class of 4-DORGs was introduced in
connection with defect tolerance schemes for nano-programmable logic
arrays~\cite{shrestha2011two}.

Finally, segment-ray graphs (SegRay) are the intersection graphs of
horizontal segments and vertical rays directed in the same direction.
SegRay graphs (and closely related graph classes) have been previously
discussed in the context of covering and hitting set problems 
(see e.g.,~\cite{Katz2005197,Chan2014112,chaplickCM13}).

In the representations defining graphs in all these classes we can assume the
$x$ and the $y$-coordinate of endpoints of any two different segments are
distinct. This property can be established by appropriate perturbations of the
segments.

The \emph{comparability graph} of a poset $P=(X,\leq_P)$ is the graph 
$(X,E)$ where for distinct $u,v\in X$ we have $uv\in E$ if and only if $u\leq_P v$ or $v\leq_P u$. 
Every bipartite graph $G=(A,B;E)$ is the comparability graph of a height-2 poset, denoted $Q_G$, where $A$ is the set of minimal elements, $B$ is the set of maximal elements, and for each $a\in A$, $b\in B$ we have $a\leq b$ in $Q_G$ if and only if $a$ and $b$ are adjacent in $G$. 
For the sake of brevity we define the \emph{dimension of a bipartite graph} $G$ to be equal to
the dimension of $Q_G$. The freedom that we may have in defining
$Q_G$, i.e., the choice of the color classes, does not affect the
dimension.
This is an easy instance of the fact that dimension is a
comparability invariant (see \cite{dim-comp-invariant}). 

% -------------------------------------------------------------
\subsection{Background on order dimension}
% -------------------------------------------------------------

Let $P = (X, \leq_P)$ be a partial order. A linear order $L = (X,
\leq_L)$ on $X$ is a \emph{linear extension} of $P$ when $x \leq_P y$
implies $x \leq_L y$.  A family ${\cal R}$ of linear extensions of $P$
is a \emph{realizer} of $P$ if $P = \bigcap_{i\in{\cal R}} L_i$, i.e., $x \leq_P
y$ if and only if $x \leq_{L} y$ for every $L \in {\cal R}$. The 
\emph{dimension} of $P$, denoted $\dim(P)$, is the minimum size of a
realizer of $P$.  This notion of dimension for partial orders was
defined by Dushnik and Miller~\cite{dm-pos-41}. The dimension of $P$
can, alternatively, be defined as the minimum $t$ such that $P$ admits
an order preserving embedding into the product order on $\RR^t$,
i.e., we can associate a $t$-vector $(x_1,\ldots,x_t)$ of
reals for each element $x\in X$ such that $x \leq_P y$ if and only
if $x_i \leq y_i$ for all $i\in\{1,\ldots,t\}$, which is denoted by
$x\leq_\dom y$.
Trotter's monograph~\cite{t-cpos-92} provides a comprehensive collection of
topics related to order dimension.

An \emph{interval order} is a partial order $P=(X,<_P)$ admitting an
interval representation, i.e., a mapping $x \to (a_x,b_x)$ from the elements
of $P$ to intervals in $\RR$ such that $x <_P y$ if and only if $b_x \leq a_y$.  The
\emph{interval dimension} of $P$, denoted $\idim(P)$, is the minimum number
$t$ such that there exist $t$ interval orders $I_i$ with $P =
\bigcap_{i=1}^{t} I_i$. Since every linear order is an interval order
$\idim(P)\leq\dim(P)$ for all $P$. If $P$ is of height two, then dimension and
interval dimension differ by at most one, i.e., $\dim(P)\leq
\idim(P)+1$,~\cite[page 47]{t-cpos-92}.

Some subclasses of grid intersection graphs are characterised by their order
dimension. For example, posets of height $2$ and dimension $2$
correspond to bipartite \emph{permutation graphs}. Bipartite permutation graphs have
an intersection representation of horizontal and vertical segments whose
endpoints lie on two parallel lines in the plane:
Drawing the first linear extension on a line and the reverse of the second
linear extension on a parallel line leads to a segment intersection
representation of the permutation graph after connecting the corresponding points
on the lines by a segment. In the bipartite case the endpoints can be arranged
on the lines such that the segments of the same color class are parallel.
Another example of a class of grid intersection graphs which is
characterised by a variant of dimension is the class of 2-DORGs.

\begin{proposition}\label{prop:2dorg+idim}
2-DORGs are exactly the bipartite graphs of interval dimension 2.
\end{proposition}
This has been shown in~\cite{STU10} using a
characterization of 2-DORGs as the complement of co-bipartite circular arc
graphs.  Below we give a simple direct proof.

To explain it we begin with a geometric version of
interval dimension:  Vectors $a,b\in \RR^d$ with $a\leq_\dom b$ define a \emph{standard box} $[a,b]= \{v : a\leq_\dom v\leq_\dom{b}\}$ in $\RR^d$.
Let $P=(X,\leq_P)$ be a poset.
A family of standard boxes $\{[a_x,b_x]\subseteq \RR^d : x\in X\}$ is a \emph{box representation} of $P$ in $\RR^d$ if it holds that $x<_P y$ if and only if $b_x\leq_\dom a_y$. 
Then the interval dimension of $P$ is the minimum $d$ for which there is a box representation of $P$ in $\RR^d$.
Note that if $P$ has height $2$ with $A=\Min(P)$ and $B=\Max(P)$, then in a box representation the lower corner $a_x$ for each $x\in A$ and the
upper corner $ b_y$ for each $y\in B$ are irrelevant, in the sense that they can uniformly be chosen as $(-c,\ldots,-c)$ respectively\ $(c,\ldots,c)$ for a large enough constant $c$.

%=============
% in einem figure environment mit caption
\calc_figscale{30}%
\begin{figure}[htb]
    \centerline{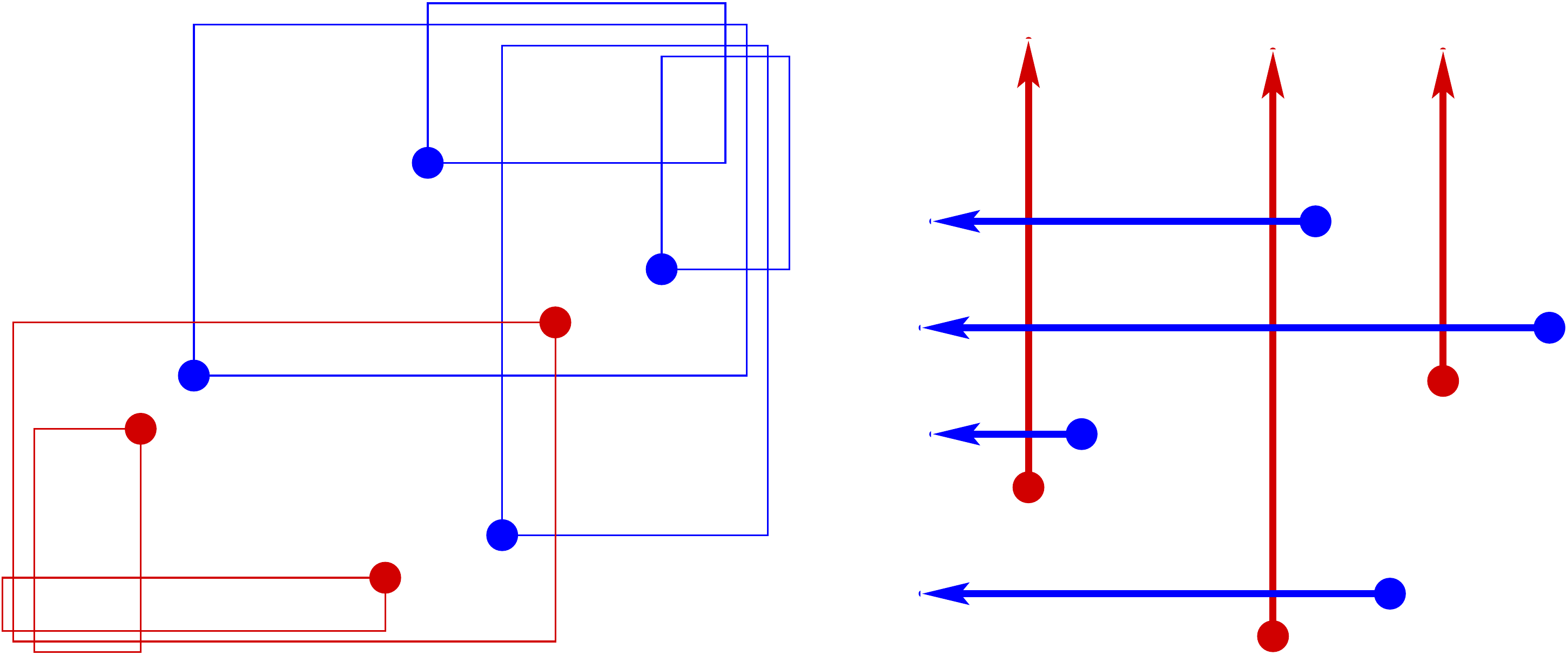}
    \caption{A representation with boxes
showing that $\idim(G)=2$ and the corresponding 2-DORG
representation.\label{fig:2dorg+idim}}
    \end{figure}%

%=============

\medskip
\ProofOf{Proposition~\ref{prop:2dorg+idim}} Let $G=(A,B;E)$ be a bipartite graph and suppose $\idim(G)=2$.
Consider a box representation $\{[a_x,b_x]\}$ of $G$ in $\RR^2$.
Clearly, for each $x\in A$, $y\in B$ we have $xy\in E$ if and only if $b_x\leq_\dom a_y$. 
To obtain a \textup{2-DORG} representation of $G$, draw upward rays starting from upper corners of boxes representing $A$, and leftward rays starting from lower corners of boxes representing $B$ (see Figure~\ref{fig:2dorg+idim}).
Then for each $x\in A$, $y\in B$ we have $b_x\leq_\dom a_y$ if and only if the rays of $x$ and $y$ intersect.

Now, the converse direction is immediate together with the observation from the previous paragraph.\qed

Since the \textup{2-DORG}s are exactly the \textup{bipartite graphs} of interval dimension
$2$, and interval dimension is bounded by dimension, we obtain
$$
\textup{2-dim BipG} \subseteq \textup{2-DORG}.
$$

% -------------------------------------------------------------
\section{Containment relations between the classes}\label{sec:relations}
% -------------------------------------------------------------
%
The diagram shown in Figure~\ref{fig:classinclusion} has 19 non-transitive
inclusions represented by the edges.  In this
section we show the inclusion between the respective classes of graphs. The
inclusion $\textup{2-dim BipG} \subseteq \textup{2-DORG}$
was already noted as a consequence of Proposition~\ref{prop:2dorg+idim}.
The next 8 inclusions follow directly from the definition of the classes:
\begin{align}
\textup{UGIG} &\subseteq \textup{GIG} &
\textup{StabGIG} &\subseteq\textup{GIG} \notag\\
\textup{2-DORG} &\subseteq \textup{3-DORG} &
\textup{3-DORG} &\subseteq \textup{4-DORG} \notag\\
\textup{3-dim GIG} &\subseteq \textup{3-dim BipG} &
\textup{3-dim BipG} &\subseteq  \textup{4-dim BipG} \notag\\
\textup{Stick} &\subseteq\textup{BipHook} &
\textup{3-dim GIG} &\subseteq \textup{GIG} \notag.
\end{align}
The following less trivial inclusions follow from
geometric modifications of the representation.
The proofs are given in the following two propositions. 
$$
\textup{BipHook}\subseteq \textup{StabGIG} \qquad\qquad \textup{2-DORG}\subseteq\textup{Stick}.
$$

\begin{figure}[htb]
    \centerline{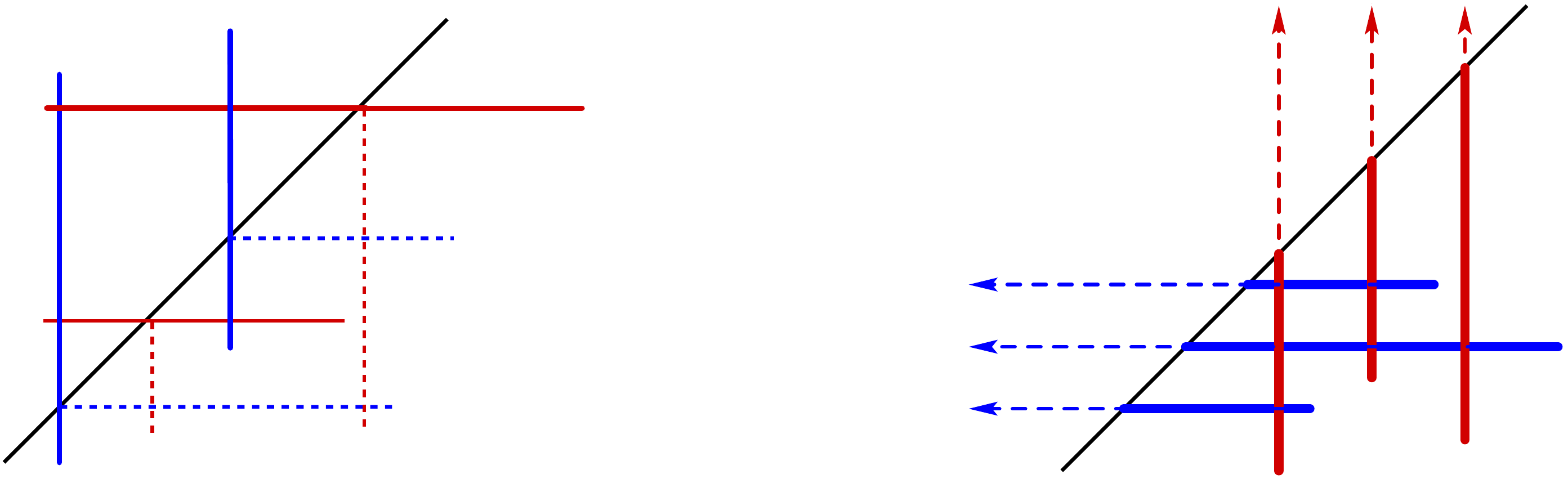}
    \caption{From a BipHook to a StabGIG  and from a 2-DORG
  to a Stick representation.\label{fig:modhookstick}}
    \end{figure}

\begin{proposition}
%$\textup{BipHook}\subseteq \textup{StabGIG}$.
Each bipartite hook graph is a stabbable GIG.
\end{proposition}
\Proof Let $G=(A,B;E)$ be a bipartite hook graph and fix a hook
representation of $G$ in which vertices of $A$ and $B$ are represented by blue and red hooks, respectively. We reflect the
horizontal part of each blue hook (dotted in
Figure~\ref{fig:modhookstick}) and the vertical part of each red hook (red dotted) at the diagonal. We claim that this
results in a StabGIG representation of the same graph. The edges are
preserved by the operation, since each intersection is witnessed by a vertical and a horizontal segment, and either both segments are reflected or none of them. On the other hand,
the transformation is an invertible linear transformation on a subset
of the segments from the region below the line to the one above,
hence no new intersection is introduced. The stabbability of the 
GIG representation comes for free. 
\qed

%%%%%%%%%%%%%%%%%%%%%%%%%%%%%%%%%%%%%%%%%%%%%%%%%%%%%%%%%%%%%%%
% in einem figure environment mit caption
\calc_figscale{30}%

%%%%%%%%%%%%%%%%%%%%%%%%%%%%%%%%%%%%%%%%%%%%%%%%%%%%%%%%%%%%%%%

\begin{proposition}\label{prop:stick2dorg}
Each 2-DORG is a Stick graph. 
\end{proposition}

\Proof
  Given a 2-DORG representation of $G$ with upward and leftward rays, we
  let $\ell$ be a line with slope 1 that is placed above all
  intersection points and endpoints of rays. 
  Removing the parts of the rays that lie in the halfplane above $\ell$ leaves a Stick representation of $G$, 
  see Figure~\ref{fig:modhookstick}.\qed

%Simple pruning
Pruning of rays also yields the following three inclusions:%13 total
$$
\textup{3-DORG} \subseteq\textup{SegRay} 
\qquad\qquad
\textup{SegRay} \subseteq \textup{GIG}
\qquad\qquad
\textup{4-DORG}\subseteq\textup{UGIG}.
$$

For the last one, consider a 4-DORG representation and a square of size $D$
that contains all intersections and endpoints of the rays.  Cutting each ray to
a segment of length $D$ leads to a UGIG representation of the same graph.
This was already observed in~\cite{STU10}.

Conversely, extending the vertical segments of a Stick representation to
vertical upward rays yields:
$$\textup{Stick} \subseteq\textup{SegRay}.$$

\begin{figure}[htb]
    \centerline{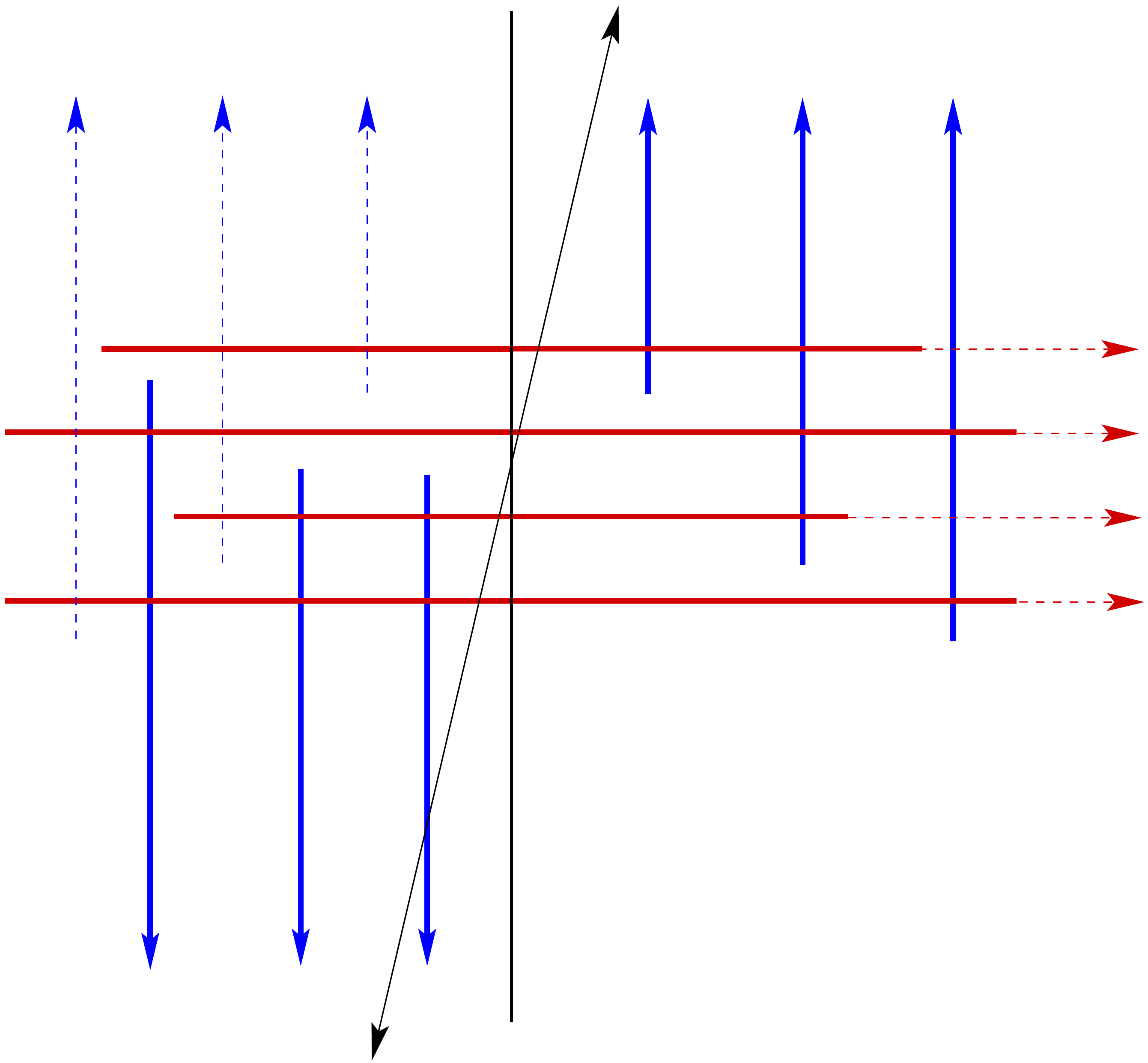}
    \caption{From a \textup{3-DORG} 
to a \textup{StabGIG} representation\label{fig:stabGIG}}
    \end{figure}

To show that every 3-DORG is a StabGIG, we use a simple geometric argument as
depicted in Figure~\ref{fig:stabGIG} and formalized in the following
proposition.
$$\textup{3-DORG}\subseteq\textup{StabGIG}.$$

%=========================
% in einem figure environment mit caption
\calc_figscale{24}%

%=========================

\begin{proposition}\label{prop:3DORGStab}
Each 3-DORG is a stabbable GIG.
%$\textup{3-DORG}\subseteq\textup{StabGIG}$.
\end{proposition}

\Proof Consider a 3-DORG representation of a graph $G$. We assume that vertical rays
point up or down while horizontal rays point right. Let $s$ be a vertical
line to the right of all the intersections.  We
prune the horizontal rays at $s$ to make them segments and then reflect the
segments at $s$, this doubles the length of the segments (see Figure~\ref{fig:stabGIG}). Now take all upward
rays and move them to the right via a reflection at $s$.  This results in an
intersection representation with vertical rays in both directions and
horizontal segments such that all rays pointing down are left of $s$ and all
rays pointing up are to the right of $s$.  Due to this property we find a
line $\ell$ of positive slope that stabs all the rays and segments of the
representation.  Pruning the rays above, respectively below their intersection with
$\ell$ yields a StabGIG representation of $G$.  \qed

A non-geometric modification of a representation gives the $16^{th}$ of the
19 non-transitive inclusions from Figure~\ref{fig:classinclusion}:

$$\textup{BipHook}\subseteq\textup{SegRay}.$$

\begin{proposition}\label{prop:HOOKinSEGRAY}
Each bipartite hook graph is a SegRay graph.
%$\textup{BipHook}\subseteq\textup{SegRay}$.
\end{proposition}

\Proof%
Consider a BipHook representation of $G=(A,B;E)$. We construct a SegRay
representation where $A$ is represented by vertical rays and $B$ by horizontal segments.
Let $a_1,\ldots,a_{|A|}$ be the order of the vertices of $A$ that we get by the centers of the hooks on the diagonal, read from bottom-left to top-right.
The $y$-coordinates of the horizontal segments and the endpoints of
the rays in our SegRay representation of $G$ will be given in the following way.

We initialize a list $R=[a_1,\dots,a_{|A|}]$ and a set
$S=B$ of \emph{active} vertices, and an empty list $Y$.
We apply one of the following steps repeatedly:
\begin{enumerate}
  \item If there is an active $a\in R$ such that $N(a)\cap S=\emptyset$,
        then remove $a$ from $R$ and append it to $Y$.
  \item \label{consecutive} If there is an active $b\in S$ such that vertices of $N(b)$
        appear consecutively in $R$, then remove $b$ from $S$ and append it to $Y$.
\end{enumerate}
Suppose that $R$ and $S$ are empty after the iteration.
Then we can construct a SegRay representation of $G$.
The endpoint of the ray representing $a_i$ receives $i$ as the $x$-coordinate and the position of $a_i$ in $Y$ as the $y$-coordinate.
The segment representing $b\in B$ also obtains the $y$-coordinate according to its position in $Y$.
Its $x$-coordinates are determined by its neighbourhood.
Now it is straight-forward to verify that this defines a SegRay representation of $G$.

%%%%%%%%%%%%%%%%%%%%%%%%%%%%%%%%%%%%%%%%%%%%%%%%%%%%%%%%%%%%%%%%%%%%%%%%%%
% in einem figure environment mit caption
\calc_figscale{40}%
\begin{figure}[htb]
    \centerline{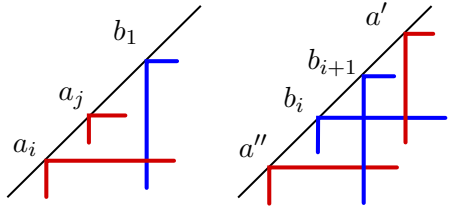}
    \caption{Left: $b_1$ is not involved in a forward pair.
Right: A forbidden configuration for a bipartite hook graph.\label{fig:HOOKtoSEGRAY}}
    \end{figure}%
%
%%%%%%%%%%%%%%%%%%%%%%%%%%%%%%%%%%%%%%%%%%%%%%%%%%%%%%%%%%%%%%%%%%%%%%%%%%

It remains to show that one of the steps can always be applied if $R$ and $S$
are nonempty.  Suppose that none of the steps can be applied.  Then, for
each $b\in S$ there are active vertices $a_i,a_k\in R\cap N(b)$ and $a_j\in
R\setminus N(b)$ with $i<j<k$.  We call $(a_j,b)$ an \emph{interesting pair}.
If the center of the hook $a_j$ lies before the center of $b$ then we call the
interesting pair \emph{forward}, and \emph{backward} otherwise. Let
$b_1,\dots,b_{|S|}$ be the order of centers of active vertices of $S$. If $b_1$
is involved in a forward interesting pair, then each certifying $a_j$ is
locked in the triangle between $b_1$ and $a_i$ and thus has no neighbour in $S$ (see Figure~\ref{fig:HOOKtoSEGRAY}), and so step 1 could be applied.
Hence every interesting pair involving
$b_1$ is a backward pair.  Symmetrically, every interesting pair involving
$b_{|s|}$ is a forward pair.  We conclude that there are active vertices
$b_i,b_{i+1}$, such that $b_i$ is involved in a forward interesting pair, and
$b_{i+1}$ in a backward one.  

Let $a'\in N(b_i)$ be the hook corresponding to an active vertex that encloses
$a\not\in N(b_i)$, i.e., $b_i<a<a'$ on the diagonal. Since $a$ is active its
hook intersects some $b_j$ with $i+1 \leq j$. Therefore, $b_i<b_{i+1}<a'$ on
the diagonal. By symmetric reasoning we also find $a''\in N(b_{i+1})$ such
that on the diagonal $a''<b_i<b_{i+1}<a'$ and  $a''b_{i+1} \in E$.
The order on the diagonal and 
the existence of edges $a''b_{i+1}$ and $b_{i}a'$ implies that the hooks of
$b_i$ and $b_{i+1}$ intersect (see Figure~\ref{fig:HOOKtoSEGRAY}).
This contradicts that $b_i$ and $b_{i+1}$ belong to the same color class of
the bipartite graph.\qed

%===================================================
\section{Dimension}\label{sec:dimension}
%===================================================

From the 19 inclusion relations between classes that have been mentioned at
the beginning of the previous section we have shown 16. The remaining three
inclusions will be shown by using order dimension in this section.  Specifically, we bound the maximal dimension of the graphs in the
relevant classes.  First, we will show that the dimension of \textup{GIG}s is
bounded.  It has previously been observed that $\idim(G) \leq 4$ when $G$ is a
GIG~\cite{ChaplickHOSU14}.  As already shown in~\cite{felsner_theorder} this
can be strengthened to $\dim(G)\leq 4$.

\begin{figure}[htb]
    \centerline{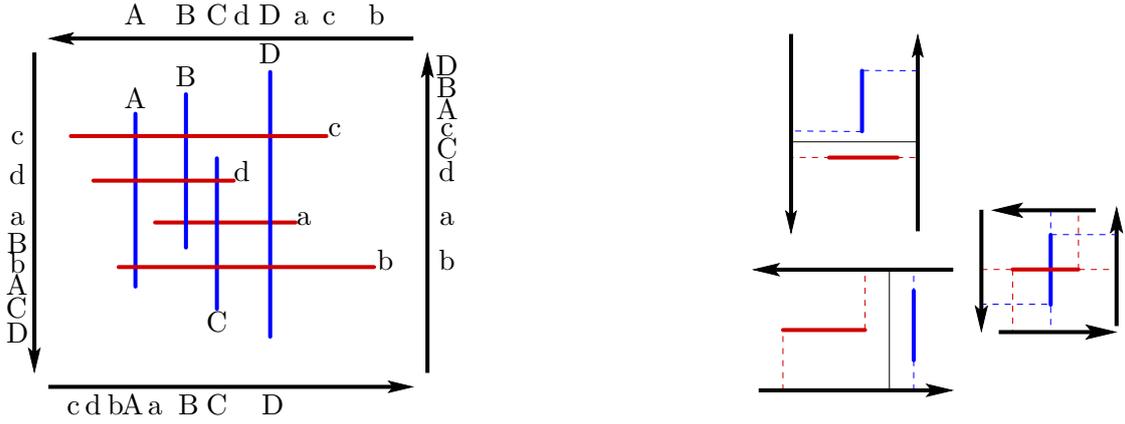}
    \caption{A realizer of a GIG and an illustration 
of the correctness.\label{fig:example4dim}}
    \end{figure}

We define four linear extensions of $G$ as depicted in
Figure~\ref{fig:example4dim}.  In each of the directions left, right, top and
bottom we consider the orthogonal projection of the segments onto a directed
horizontal or vertical line. In each such projection every segment corresponds to
one interval (or point) per line.  We choose a point from each interval on the
line by the following rule.  For minimal elements we take the minimal point in
the interval in the direction of the line, for maximal elements we choose the
maximal one.  We denote those total orders according to the direction of their
oriented line by $L_\leftarrow, L_\rightarrow, L_\uparrow, L_\downarrow$.

\begin{proposition}\label{prop:GIGdim4}
For every GIG $G$, $\{L_\leftarrow, L_\rightarrow, L_\uparrow, L_\downarrow\}$ is a realizer of $G$. 
% Let $G$ be a GIG. 
% Then $\{L_\leftarrow, L_\rightarrow, L_\uparrow, L_\downarrow\}$ is a realizer of $G$.
 Hence \mbox{$\dim(G)\leq 4$}.
\end{proposition}
%
%%%%%%%%%%%%%%%%%%%%%%%%%%%%%%%%%%%%%%%%%%%%%%%%%%%%%%%%%%%%%%%%%%%%%%%%%%
% in einem figure environment mit caption
\calc_figscale{40}%

%%%%%%%%%%%%%%%%%%%%%%%%%%%%%%%%%%%%%%%%%%%%%%%%%%%%%%%%%%%%%%%%%%%%%%%%%%
\Proof
For two intersecting segments the minimum always lies before the maximum, see
Figure~\ref{fig:example4dim}.  It remains to check that every incomparable
pair $(s_1,s_2)$ is reversed in the realizer.  Every disjoint pair of segments
is separated by a horizontal or vertical line.  The separated vertices
appear in different order in the two directions orthogonal to this line, see
Figure~\ref{fig:example4dim}.\qed

It is known that a bipartite graph is a bipartite permutation graph if and only if the
dimension of the poset is at most $2$. Thus, by Proposition~\ref{prop:GIGdim4}, the maximal dimension of the graphs
in the other classes that we consider must be $3$ or $4$.  In the following we show that BipHooks
and 3-DORGs have dimension at most $3$. For these results we note that graphs with a special SegRay representation have dimension at most~$3$ and that the interval dimension of SegRay graphs is bounded by~$3$. This latter result was shown
previously by a different argument in~\cite{ChaplickHOSU14}.

\begin{lemma}\label{lem:segrayprop}
%Let $G$ be a SegRay graph with a SegRay representation satisfying the following:
%whenever two horizontal segments are such that the $x$-projection of one is included in the other one, then the smaller segment lies below the bigger one. Then $\dim(G)\leq 3$.
For a SegRay graph $G$, \mbox{$\dim(G) \leq 3$} when $G$ has a SegRay representation satisfying the following:
whenever two horizontal segments are such that the $x$-projection of one is included in the other one, then the smaller segment lies below the bigger one.
\end{lemma}
\Proof
Consider such a SegRay representation of $G$ with horizontal segments as
maximal and downward rays as minimal elements of $Q_G$.
The linear extensions $L_\rightarrow$, $L_\leftarrow$, and $L_\downarrow$ defined for
Proposition~\ref{prop:GIGdim4} form a realizer of $Q_G$.
\qed

\begin{corollary}
 %Let $G$ be a 3-DORG. Then $\dim(G)\leq 3$.
 For every 3-DORG $G$, $\dim(G) \leq 3$. 
\end{corollary}
\Proof Consider a 3-DORG representation of $G$ where the horizontal rays use
two directions.  We cut the horizontal rays so that they have the same length
$D$. When $D$ is large enough, this yields a SegRay representation of the same graph. Note that such a representation has no nested segments. Thus, Lemma 1 implies $dim(Q_G) \leq 3$. \qed

\begin{proposition}\label{prop:segrayidim}
% Let $G$ be a SegRay graph. Then $\idim(G)\leq 3$.
 For every SegRay graph $G$, $\idim(G) \leq 3$. 
\end{proposition}
\Proof Suppose that the rays correspond to minimal elements of $Q_G$.
By Lemma~\ref{lem:segrayprop} the linear extensions $L_\rightarrow$, $L_\leftarrow$ and $L_\downarrow$ reverse all incomparable pairs except some that consist of two maximal elements.  We convert these linear extensions to interval orders and extend the intervals (originally points) of maximal elements in
$L_\rightarrow$ far to the right to make them intersect. In this way we obtain three interval orders whose intersection gives rise to $Q_G$.
\qed

\begin{proposition}
%Let $G$ be a bipartite hook graph. Then $\dim(G)\leq 3$.
For every bipartite hook graph $G$, $\dim(G) \leq 3$. 
\end{proposition}
\Proof Let $A$ and $B$ be the color classes of $G$.
We construct the graph $G'$ by adding private neighbours to vertices of $B$.
Then $G'$ is also a BipHook graph as we can easily add hooks intersecting a single hook in a representation of $G$.
By Proposition~\ref{prop:HOOKinSEGRAY} we know that $G'$ has a SegRay representation $R$ with downward rays representing $A$.
By construction, each horizontal segment in $R$ must have its private ray intersecting it.
Thus $R$ satisfies the property of Lemma~\ref{lem:segrayprop} and $\dim(Q_{G'})\leq 3$.
Since $Q_G$ is an induced subposet of $Q_{G'}$ we conclude $\dim(Q_G)\leq 3$.\qed

Since $\textup{Stick}\subset \textup{BipHook}$ we know that $\dim(Q_G)\leq 3$ if
$G$ is a Stick graph.  However, a nicer realizer for a Stick graph is obtained
by Proposition~\ref{prop:GIGdim4}, since $L_\leftarrow$ and $L_\downarrow$
coincide in a Stick representation.

In Section~\ref{sec:separating} we will show that these bounds are tight.

%==================
\section{Vertex-Edge Incidence Posets}\label{sec:incidence}
%==================

We proceed by investigating the relations between the classes of GIGs and 
incidence posets of graphs.

For a graph $G$, $P_G$ denotes the vertex-edge incidence poset of $G$, and the
comparability graph of $P_G$ is the graph obtained by subdividing each edge of
$G$ once.  The vertex-edge incidence posets of dimension $3$ are characterised by
Schnyder's Theorem.

\begin{theorem}[\cite{schnyder1989planar}]
A graph $G$ is planar if and only if $\dim(P_G)\leq 3$.
\end{theorem}
Even though some GIGs have poset dimension $4$, we will see that the
vertex-edge incidence posets with a GIG representation are precisely the vertex-edge incidence posets of
planar graphs.

A \emph{weak bar-visibility representation} of a graph is 
a drawing that represents the vertices as horizontal segments and the edges 
as vertical segments (sight lines) touching its adjacent vertices. 

\begin{theorem}[\cite{tamassia1986unified,wismath1985characterizing}]
\label{thm:barvisibility}
A graph $G$ has a weak bar-visibility representation if and only if $G$ is planar. 
\end{theorem}

A weak bar-visibility representation of $G$ gives a GIG representation of
$P_G$.  On the other hand, a GIG representation of $P_G$ can be transformed
into a weak bar-visibility representation of $G$.  In particular, since the segments
representing edges of $G$ intersect two segments representing incident vertices, they can be shortened
until their intersections become contacts.  Hence $P_G$ is a GIG if and only if $G$ is
planar.  We next show the stronger result, that there is a $\mbox{StabGIG}$
representation of $P_G$ for every planar graph $G$.
\begin{proposition}\label{prop:planarstabbable}
A graph $G$ is planar if and only if $P_G$ is a stabbable GIG.
\end{proposition}

We use the following definitions.  A \emph{generic floorplan} is a partition
of a rectangle into a finite set of interiorly disjoint rectangles that have
no point where four rectangles meet.  Two floorplans are \emph{weakly
  equivalent} if there exist bijections $\Phi_H$ between the horizontal
segments and $\Phi_V$ between the vertical segments, such that a segment $s$
has an endpoint on $t$ in $F$ if and only if $\Phi(s)$ has an endpoint on $\Phi(t)$. A
floorplan $F$ \emph{covers} a set of points $P$ if and only if every segment contains
exactly one point of $P$ and no point is contained in two segments. The
following theorem has been conjectured by Ackerman, Barequet and
Pinter~\cite{ackerman2006number}, who have also shown it for the special case
of \emph{separable} permutations.  It has been shown by
Felsner~\cite{felsnerexploiting} for general permutations.

\begin{theorem}[\cite{felsnerexploiting}]\label{thm:floorplan}
 Let $P$ be a set of $n$ points in the plane, such that no two points have the
 same $x$- or $y$-coordinate and $F$ a generic floorplan with $n$ segments.
 Then there exists a floorplan $F'$, such that $F$ and $F'$ are weakly
 equivalent and $F'$ covers $P$. 
\end{theorem}

\ProofOf{Proposition~\ref{prop:planarstabbable}} Consider a weak bar-visibility
representation of $G$. The lowest and highest horizontal segments $h_b$ and
$h_t$ can be extended, such that their left as well as their right endpoints
can be connected by new vertical segments $v_l$ and $v_r$. The segments $h_b$,
$h_t$, $v_l$ and $v_r$ are the boundary of a rectangle. Extending every
horizontal segment until its left and right endpoints touch vertical segments
leads to a floorplan $F$. By Theorem~\ref{thm:floorplan} there exists an
equivalent floorplan $F'$ that covers a pointset $P$ consisting of $n$ points
on the diagonal of the the big rectangle with positive slope. Shortening the horizontal segments and extending the vertical segments of $F'$ by $\epsilon>0$
on each end leads to a GIG
representation of $P_G$ that can be stabbed by the line through the diagonal.

On the other hand, every GIG representation of $P_G$ leads to a weak
bar-visibility representation, and hence $G$ is planar.\qed

We will now show that $P_G$ is in the classes of Stick and 
bipartite hook graphs if and only if $G$ is outerplanar.

\begin{proposition}\label{prop:stickouterplanar}
$P_G$ is a Stick graph if and only if $G$ is outerplanar.
\end{proposition}

\Proof
Outerplanar graphs have been characterized by linear orderings of their
vertices by Felsner and Trotter~\cite{felsner2005posets}: A graph $G=(V,E)$ is
outerplanar if and only if there exist linear orders $L_1,L_2,L_3$ of the vertices with
$L_2=\overleftarrow{L_1}$, i.e., $L_2$ is the reverse of $L_1$, such that for
each edge $vw\in E$ and each vertex $u\not\in \{v,w\}$ there is $i\in\{1,2,3\}$, such
that $u>v$ and $u>w$ in $L_i$. 

Consider a Stick representation of $P_G$ where the elements of $V$ correspond to vertical sticks. Restricting the linear
extensions $L_1=L_\leftarrow$, $L_2=L_\rightarrow$, and $L_3=L_\uparrow$
(cf.~the proof of Proposition~\ref{prop:GIGdim4}) obtained from a Stick
representation of $P_G$ to the elements of $V$ yields linear orders satisfying the property above. Thus $G$ is outerplanar.

\begin{figure}[htb]
    \centerline{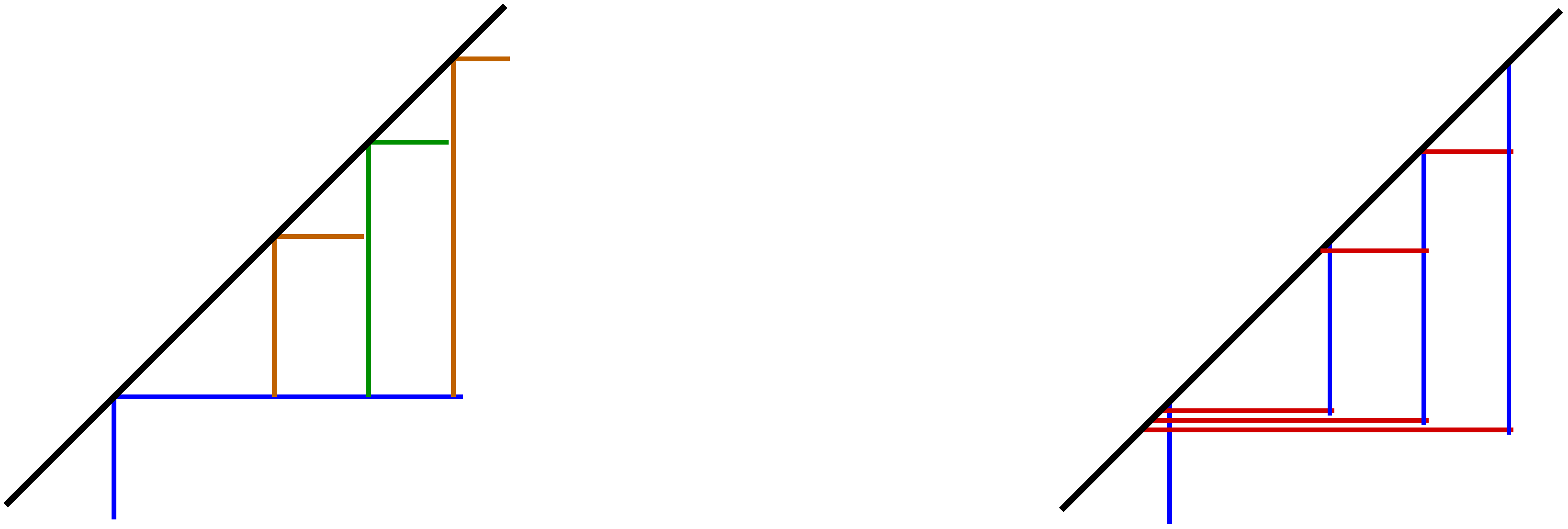}
    \caption{A hook contact representation of $G$
transformed into a Stick representation of $P_G$.\label{fig:stickP_G}}
    \end{figure}

For the backward direction let $G$ be an outerplanar graph.  In~\cite{maxpointtolerance} it is shown
that the class of {\em hook contact} graphs (each intersection of hooks is
also an endpoint of a hook) is exactly the class of outerplanar graphs. Given a
hook contact representation of $G$ we construct a Stick representation of
$P_G$. To this end we consider each hook as two sticks, a vertical one for the
vertices and a horizontal one as a placeholder for the edges.  For each
contact of the horizontal part of a hook $v$ we place an additional horizontal
stick slightly below the center of $v$. The $k$-th contact of a hook with the
horizontal part is realized by the $k$-th highest edge that is added in the
placeholder as shown in Figure~\ref{fig:stickP_G}.\qed

%%%%%%%%%%%%%%%%%%%%%%%%%%%%%%%%%%%%%%%%%%%%%%%%%%%%%%%%
% in einem figure environment mit caption
\calc_figscale{30}%

%%%%%%%%%%%%%%%%%%%%%%%%%%%%%%%%%%%%%%%%%%%%%%%%%%%%%%%%

We continue by providing some characterizations of outerplanar graphs according to GIG representations of vertex-edge incidence posets.

A \emph{weak semibar-visibility representation} of a graph is a drawing 
that represents the vertices as vertical segments with lower end at the horizontal line
$y=0$, and the edges as horizontal segments touching the two vertical segments that represent incident vertices.
 
\begin{lemma}[\cite{cobos1996visibility}]\label{lem:semibar}
  A graph $G$ is outerplanar if and only if $G$ has a weak semibar-visibility representation.
\end{lemma}

The construction used in the previous proposition directly produces a weak
semibar-visibility representation of an outerplanar graph.  Just extend all
vertical segments upwards until they hit a common horizontal line $\ell$ and
reflect the plane at $\ell$, now $\ell$ can play the role of the $x$-axis for
the weak semibar-visibility representation.

\begin{proposition}\label{prop:segrayouterplanar}
A graph $G$ is outerplanar if and only if the graph $P_G$ has a SegRay representation where the vertices of $G$ are represented as rays.
\end{proposition}
\Proof Cutting the rays of a SegRay representation with rays pointing
downwards somewhere below all horizontal segments leads to a weak semibar-visibility
representation of $G$ and vice versa. Thus, Lemma~\ref{lem:semibar} gives the
result.\qed
\begin{proposition}\label{prop:hookouterplanar}
A graph $G$ is outerplanar if and only if the graph $P_G$ has a hook representation.
\end{proposition}
\Proof
If $G$ is outerplanar then $P_G$ has a hook representation by
Proposition~\ref{prop:stickouterplanar}.
On the other hand, assume that $P_G$ has a hook representation for a graph $G$.
According to Proposition~\ref{prop:HOOKinSEGRAY} we construct a SegRay
representation with vertices as rays and edges as segments. This representation
shows that $G$ is outerplanar by Proposition~\ref{prop:segrayouterplanar}. 
\qed
\begin{proposition}\label{prop:segrayouterplanar2}
If $G$ is outerplanar, then the graph $P_G$ has a SegRay representation where the
vertices of $G$ are represented as segments.
\end{proposition}
\Proof Consider a hook representation $R'$ of $P_G$.  According to the proof
Proposition~\ref{prop:HOOKinSEGRAY} we can transform $P_G$ into a SegRay
representation with a free choice of the colorclass that is represented by
rays. Choosing the subdivision vertices as rays leads to the required
representation.
\qed
%

%%%%%%%%%%%%%%%%%%%%%%%%%%%%%%%%%%%%%%%%%%%%%%%%%%%%%%%%
% in einem figure environment mit caption
\calc_figscale{40}%
\begin{figure}[htb]
    \centerline{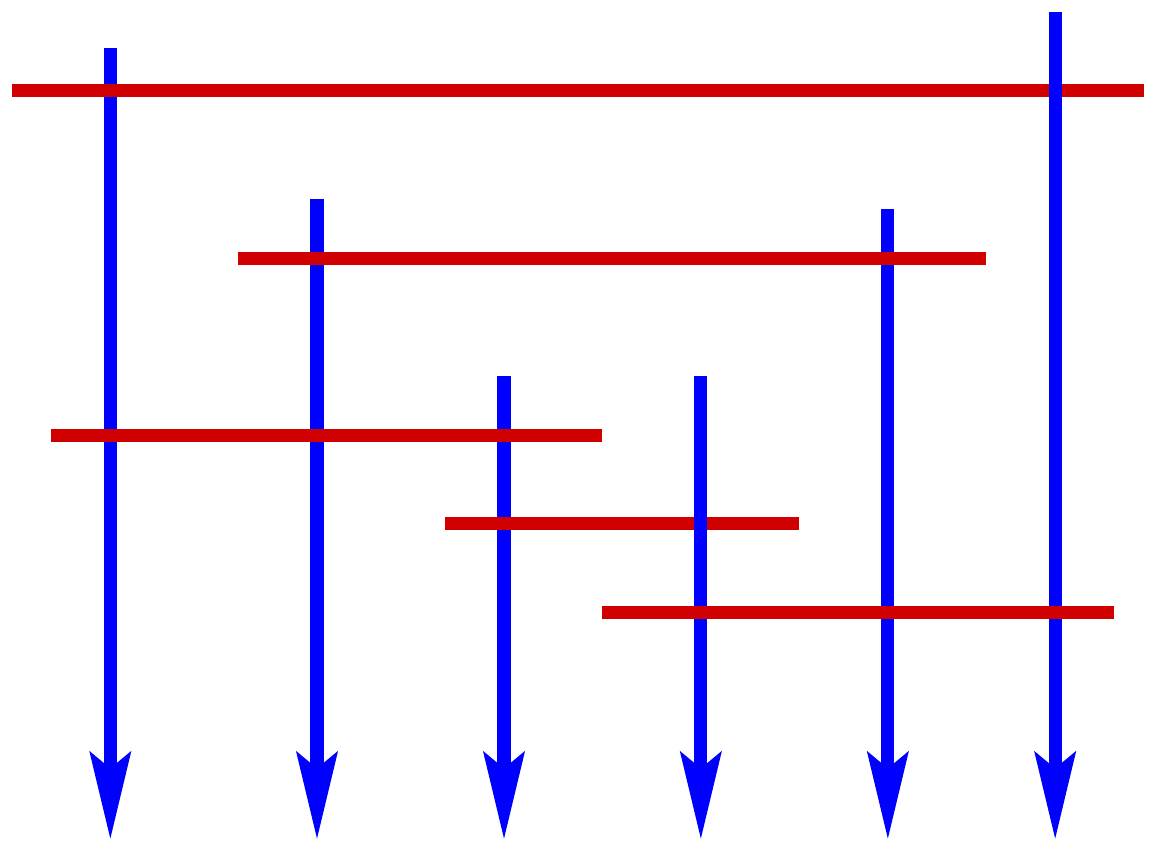}
    \caption{A SegRay representation of $P_{K_{2,3}}$.\label{fig:segray_k23}}
    \end{figure}%

%%%%%%%%%%%%%%%%%%%%%%%%%%%%%%%%%%%%%%%%%%%%%%%%%%%%%%%%
In contrast to Proposition~\ref{prop:segrayouterplanar}, the backward direction of Proposition~\ref{prop:segrayouterplanar2} does not hold:
Figure~\ref{fig:segray_k23} shows a SegRay representation of $P_{K_{2,3}}$ with vertices being represented as horizontal segments, but $K_{2,3}$ is not outerplanar.
Together with Proposition~\ref{prop:segrayouterplanar} this also shows that the class of SegRay graphs is not symmetric in its color classes.

In the following we construct a UGIG representation of $P_G$ for an
outerplanar graph $G$.
\begin{proposition}\label{prop:outerplanarUGIG}
If $G$ is outerplanar then $P_G$ is a UGIG.
\end{proposition}

\Proof We construct a UGIG representation of $P_G$ for a maximal outerplanar
graph $G=(V,E)$ with outer-face cycle $v_0,\dots,v_n$. The vertices of $V$ are
drawn as vertical segments.  Starting from $v_0$ we iteratively draw the
vertices of breadth-first-search layers (BFS-layers). Each BFS-layer has 
a natural order inherited from the order on the outer-face, i.e., the 
increasing order of indices. When the $i$-th layer $L_i$ has been
drawn the following invariants hold:
\begin{enumerate}
  \item Segments for all vertices and edges of $G[L_0,\dots,L_{i-1}]$, 
   all vertices of $L_i$, and all edges connecting vertices of $L_{i-1}$ to vertices of $L_i$ have been placed.
  \item The upper endpoints of the segments representing vertices in
  $L_i$ lie on a strict monotonically decreasing curve $C_i$. Their order on $C_i$ agrees with the order of the corresponding vertices in $L_i$.
  Their $x$-coordinates differ by at most one.
  \item No segment intersects the region above $C_i$.
 \end{enumerate}
%%%%%%%%%%%%%%%%%%%%%%%%%%%%%%%%%%%%%%%%%%%%%%%%%%%%%%%%
% in einem figure environment mit caption
\calc_figscale{60}%
\begin{figure}[htb]
    \centerline{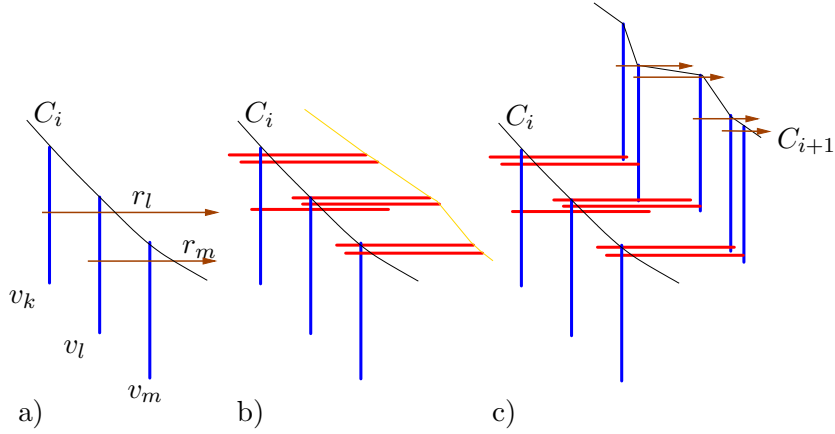}
    \caption{One step in the construction of a UGIG representation
of $P_G$: a)The situation before the step. b) The edges between layer $L_i$ and
$L_{i+1}$ and within layer $L_i$ are added. c) The vertices of layer $L_{i+1}$
are added.\label{fig:outerUGIG}}
    \end{figure}%

%%%%%%%%%%%%%%%%%%%%%%%%%%%%%%%%%%%%%%%%%%%%%%%%%%%%%%%%

We start the construction with the vertical segment corresponding to $v_0$.
The curve $C_0$ is chosen as a line with negative slope that intersects the
upper endpoint of $v_0$.

We start the $(i+1)$-th step by adding segments for the edges within vertices of layer
$L_i$.  Afterwards we add the segments for edges between vertices in layer
$L_i$ and $L_{i+1}$ and the segments for the vertices of layer $L_{i+1}$. The
construction is indicated in Figure~\ref{fig:outerUGIG}.

First we draw unit segments for the edges within layer $L_i$.  Since the graph
is outerplanar such edges only occur between consecutive vertices of the
layer. For a vertex $v_k$ of $L_i$ which is not the first vertex of $L_i$ we
define a horizontal ray $r_k$ whose start is on the segment of
the predecessor of $v_k$ on this layer such that the only additional
intersection of $r_k$ is with the segment of $v_k$. 
The initial unit segment of ray $r_k$
can be used for the edge between $v_k$ and its predecessor.

All segments that will represent edges between layer $L_i$ and $L_{i+1}$ are
placed as horizontal segments that intersect the segment of the incident
vertex $v_k\in L_i$ above the ray~$r_k$. We draw these edge-segments such that
the endpoints lie on a monotonically decreasing curve~$C$ and the order of
these endpoints on $C$ corresponds to the order of their incident vertices
in~$L_{i+1}$.

Now the right endpoints of the edges between the two layers lie on the monotone
curve $C$ and no segment intersects the region above this curve.  
Due to properties of the BFS for outerplanar graphs, each vertex of 
layer $L_{i+1}$ is incident to one or two edges whose segments end on $C$
and if there are two then they are consecutive on $C$.
We place the unit segments of vertices of $L_{i+1}$, such that 
their lower endpoint is on the lower segment of an incident edge with 
the $x$-coordinate such that they realize the required intersections.

With this construction the invariants are satisfied.
\qed

\begin{figure}[htb]
    \centerline{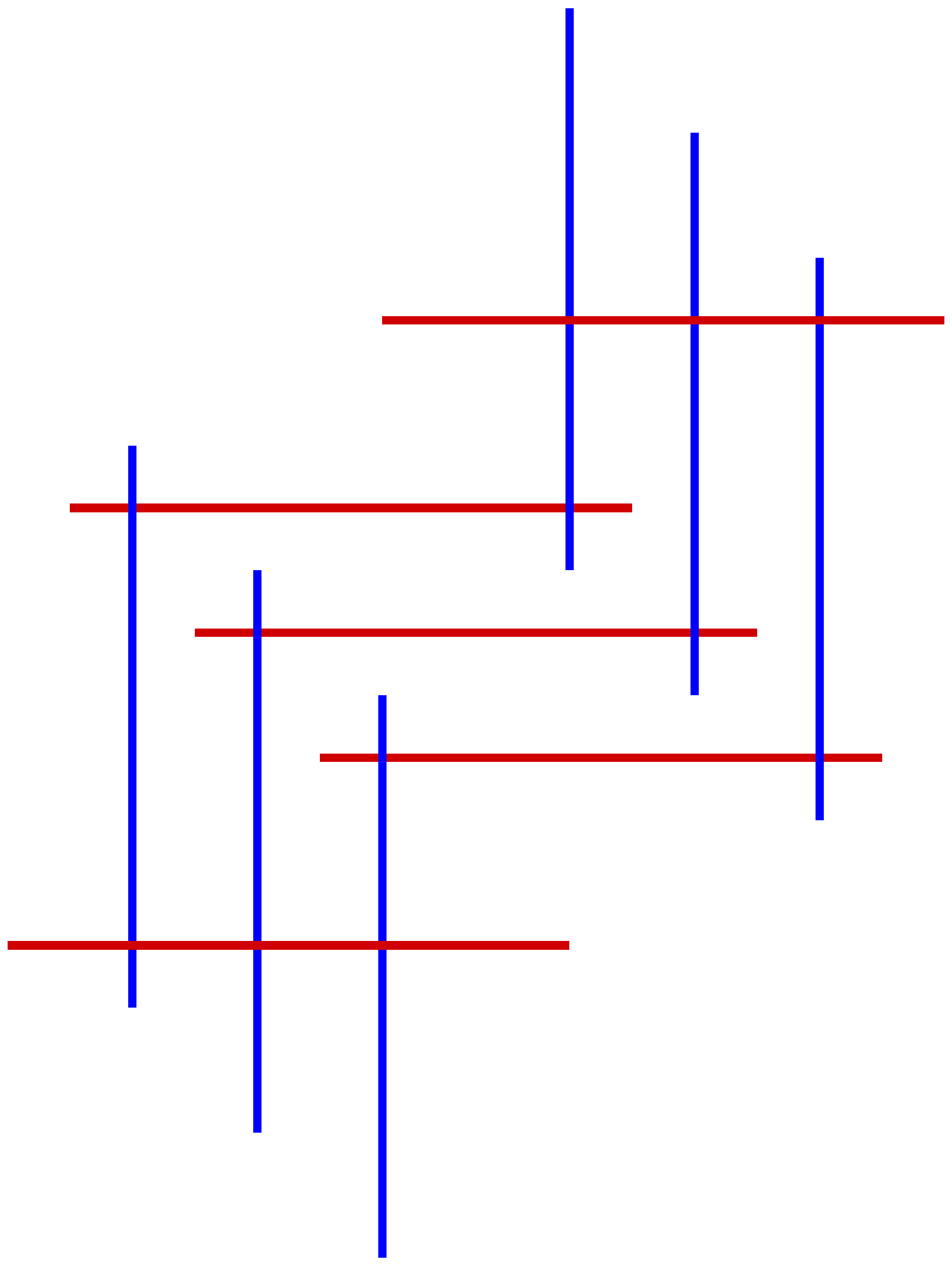}
    \caption{A UGIG representation of $P_{K_{2,3}}$.\label{fig:K23UGIG}}
    \end{figure}

There are graphs $G$ where $P_G$ is a UGIG and $G$ is not outerplanar, for
example $G=K_{2,3}$ as shown in Figure~\ref{fig:K23UGIG}.
On the other hand there exist planar graphs $G$, such that $P_{G}$ is not a
UGIG as the following proposition shows.

\begin{proposition}\label{prop:notUGIG_edgevertex}
$P_{K_4}$ is not a UGIG.
\end{proposition}

\Proof Suppose to the contrary that $P_{K_4}$ has a UGIG representation with vertices as vertical segments. 
By contracting vertical segments to points one can obtain a planar embedding of $K_4$ from such a representation.
As $K_4$ is not outerplanar, there is a vertex $v$ that is not incident to the outer face in this embedding.
For the initial UGIG representation this means that $v$ is represented by a vertical segment which is enclosed by segments representing vertices and edges of $K_4-\{v\}$.
Notice that these segments represent a 6-cycle of $P_{K_4}$.
However, the largest vertical distance between any
pair of horizontal segments in this cycle is less than~1. Thus, there is not
enough space for the vertical segment of $v$, contradiction. \qed

%%%%%%%%%%%%%%%%%%%%%%%%%%%%%%%%%%%%%%%%%%%%%
% in einem figure environment mit caption
\calc_figscale{25}%

%%%%%%%%%%%%%%%%%%%%%%%%%%%%%%%%%%%%%%%%%%%%%

%--------------------------------------------------
\section{Separating examples}\label{sec:separating}
%--------------------------------------------------

In this section we will give examples of graphs that separate the graph classes
in Figure~\ref{fig:classinclusion}. For this purpose we will show that the classes
we have observed to be at most $4$-dimensional indeed contain $4$-dimensional graphs.
This is done in Subsection~\ref{subsec:sepDimension} using standard examples and vertex-face incidence posets of outerplanar graphs. 
The remaining graph classes will
be separated using explicit constructions in
Subsection~\ref{subsec:sepConstructions} and Subsection~\ref{subsec:stab}.

Using the observations of Section~\ref{sec:incidence} about vertex-edge incidence posets we can immediately separate the following graph classes.
\begin{align*}
 \textup{StabGIG} &\not\subset \textup{BipHook} &
 \textup{StabGIG} &\not\subset \textup{3-DORG} \notag\\
 \textup{SegRay}  &\not\subset \textup{4-DORG}  & 
 \textup{Stick}   &\not\subset \textup{2-DORG}.
\end{align*}

In~\cite{STU10} it is shown that the graph $C_{14}$ (cycle on $14$ vertices) is not a $4$-DORG, and in particular is not a $3$- or $2$-DORG.
In other words, $P_{C_7}$ is not a $4$-DORG.
Since $C_7$ is outerplanar, by the propositions of the previous section we know that $P_{C_7}$ is a SegRay, a StabGIG and a Stick graph.
This shows the three seperations involving DORGs.
For the first one let $G$ be a planar graph that is not outerplanar.
Then $P_G$ is a StabGIG (Proposition~\ref{prop:planarstabbable}) but not a BipHook graph (Proposition~\ref{prop:hookouterplanar}).

%-------------------------------------------------
\subsection{4-Dimensional Graphs \label{subsec:sepDimension}}
%-------------------------------------------------

%%%%%%%%%%%%%%%%%%%%%%%%%%%%%%%%%%%%%%%%%%%%%%%%%%%%%%%%%%%%%%%%%%%%%%%%%%
% in einem figure environment mit caption
\calc_figscale{35}%
\begin{figure}[htb]
    \centerline{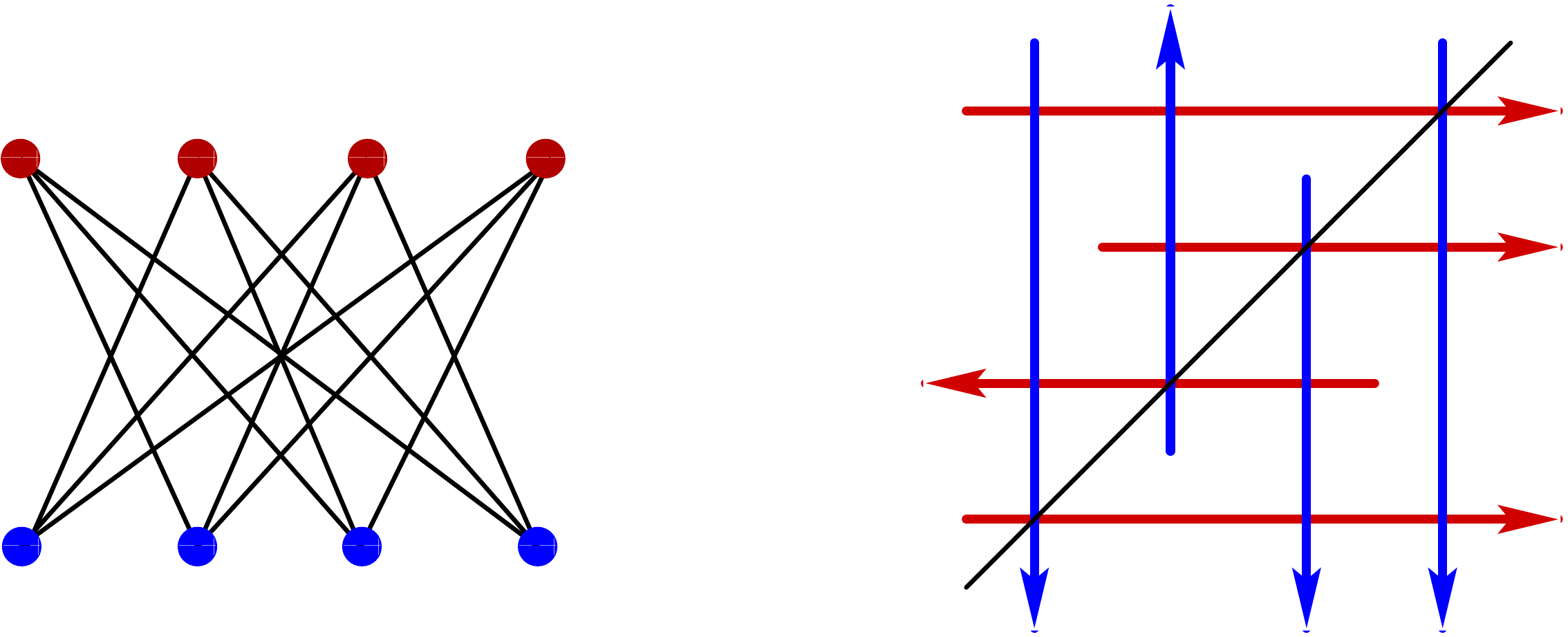}
    \caption{The poset $S_4$ and a stabbable 4-DORG representation of it.\label{fig:stabS_4}}
    \end{figure}%

%%%%%%%%%%%%%%%%%%%%%%%%%%%%%%%%%%%%%%%%%%%%%%%%%%%%%%%%%%%%%%%%%%%%%%%%%%

First of all, some graph classes are already separated by their
maximal dimension.  The \emph{standard example} $S_n$ of an $n$-dimensional
poset, cf.~\cite{t-cpos-92}, is the poset on $n$ minimal elements $a_1,\ldots,a_n$ and $n$ maximal elements $b_1,\ldots,b_n$, such that $a_i< b_j$ in $S_n$ if and only if $i\neq j$.
To separate most of the $4$-dimensional
classes from the $3$-dimensional ones, the standard example $S_4$ is sufficient.
As shown in Figure~\ref{fig:stabS_4} it has as a stabbable 4-DORG representation.
From this it follows that:
\begin{align*}
 \textup{StabGIG} &\not\subset \textup{BipHook} &
\textup{StabGIG} &\not\subset \textup{3-DORG} \notag\\
\textup{4-DORG} &\not\subset \textup{3-DORG} &
\textup{StabGIG} &\not\subset \textup{3-dim GIG}\notag
\end{align*}
Since the interval dimension of $S_n$ is $n$ 
we get the following relations from 
Proposition~\ref{prop:segrayidim}.
\begin{align*}
\textup{StabGIG} &\not\subset \textup{SegRay} &
\textup{4-DORG} &\not\subset \textup{SegRay}
\end{align*}
We will now show that the vertex-face incidence poset of an outerplanar graph has a
SegRay representation.  In~\cite{felsner2011order} it has been shown that
there are outerplanar maps with a vertex-face incidence poset of dimension $4$.  Together
with Proposition~\ref{prop:segrayoutermap} below this shows that there are
SegRay graphs of dimension $4$.
We obtain
\begin{align*}
\textup{SegRay} &\not\subset \textup{3-dim GIG}.
\end{align*}

\vskip-7mm

\begin{proposition}\label{prop:segrayoutermap}
If $G$ is an outerplanar map then the vertex-face incidence poset of $G$ is a SegRay graph.
\end{proposition}

Let $G$ be a graph with a fixed outerplanar embedding.
First we argue that we may assume that $G$ is $2$-connected.
If $G$ is not connected then we can add a single
edge between two components without changing the vertex-face poset.
Now consider adding an edge between two neighbours of a cut vertex on the
outer face cycle, i.e., two vertices of distance $2$ on this cycle.
This adds a new face to the vertex-face-poset, but keeps the old vertex-face-poset
as an induced subposet.
Therefore, we may assume that $G$ is $2$-connected.

\begin{figure}[h]
 \centering
 \includegraphics{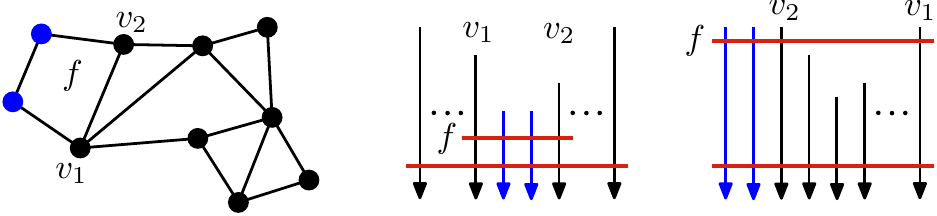}
 \caption{Illustration for the induction step in Proposition~\ref{prop:segrayoutermap}}
 \label{fig:vertex-face-segray}
\end{figure}

By induction on the number of bounded faces we show that $G$ has a SegRay representation in which the cyclic order of the vertices on the outer face agrees with the left-right order (cyclically) of rays representing these vertices.
If $G$ has one bounded face then the claim is straight-forward.
If $G$ has more bounded faces then consider the dual graph of $G$ without the outer face, which is a tree.
Let $f$ be a face that corresponds to a leaf of that tree.
Define $G'$ to be the plane graph obtained by removing $f$ and incident degree-$2$ vertices from $G$.
Then exactly two vertices $v_1,v_2$ of $f$ are still in $G'$, and they are adjacent via an edge at the outer face of $G'$.
Note that $G'$ is $2$-connected.
Applying induction on $G'$ we obtain a SegRay representation in which the two rays representing $v_1$ and $v_2$ are either consecutive, or left- and rightmost ray.

In the first case we insert rays for the removed vertices between $v_1$ and $v_2$ with endpoints being below all other horizontal segments.
Then a segment representing $f$ can easily be added to obtain a SegRay representation with the required properties of $G$, see the middle of Figure~\ref{fig:vertex-face-segray}.

If the rays of $v_1$ and $v_2$ are the left- and rightmost ones, then observe that the endpoints of both rays can be extended upwards to be above all other endpoints.
We can insert the new rays to the left of all the other rays and the segment for $f$ as indicated in Figure~\ref{fig:vertex-face-segray} on the right.
This concludes the proof.\qed

Propositions~\ref{prop:segrayoutermap} and~\ref{prop:segrayidim} also give the
following interesting result about vertex-face incidence posets of outerplanar maps
which complements the fact that they can have dimension $4$ \cite{felsner2011order}.

\begin{corollary}\label{cor:idim-opl}
 The interval dimension of a vertex-face incidence poset of an outerplanar 
 map is bounded by $3$.
\end{corollary}

We have separated all the graph classes which involve dimension except for the two classes of $3$-dimensional GIGs and stabbable GIGs.
As indicated in Figure~\ref{fig:classinclusion} it remains open whether $3$-dim GIG is a subclass of StabGIG or not.
More comments on this can be found at the end of Subsection~\ref{subsec:stab}.

\goodbreak

\subsection{Constructions \label{subsec:sepConstructions}}

In this subsection we give explicit constructions for the remaining separations of classes not involving StabGIG.

In the introduction we mentioned that every 2-dimensional order of height $2$,
i.e., every bipartite permutation graph, is a GIG.
We show now that this does not hold for 3-dimensional orders of height $2$.
\begin{proposition}\label{prop:threedimnonpseudoseg}
There is a 3-dimensional bipartite graph that is not a GIG.
\end{proposition}
%
%%%%%%%%%%%%%%%%%%%%%%%%%%%%%%%%%%%%%%%%%%%%%%%%%%%%%%%%
% in einem figure environment mit caption
\calc_figscale{25}%
\begin{figure}[htb]
    \centerline{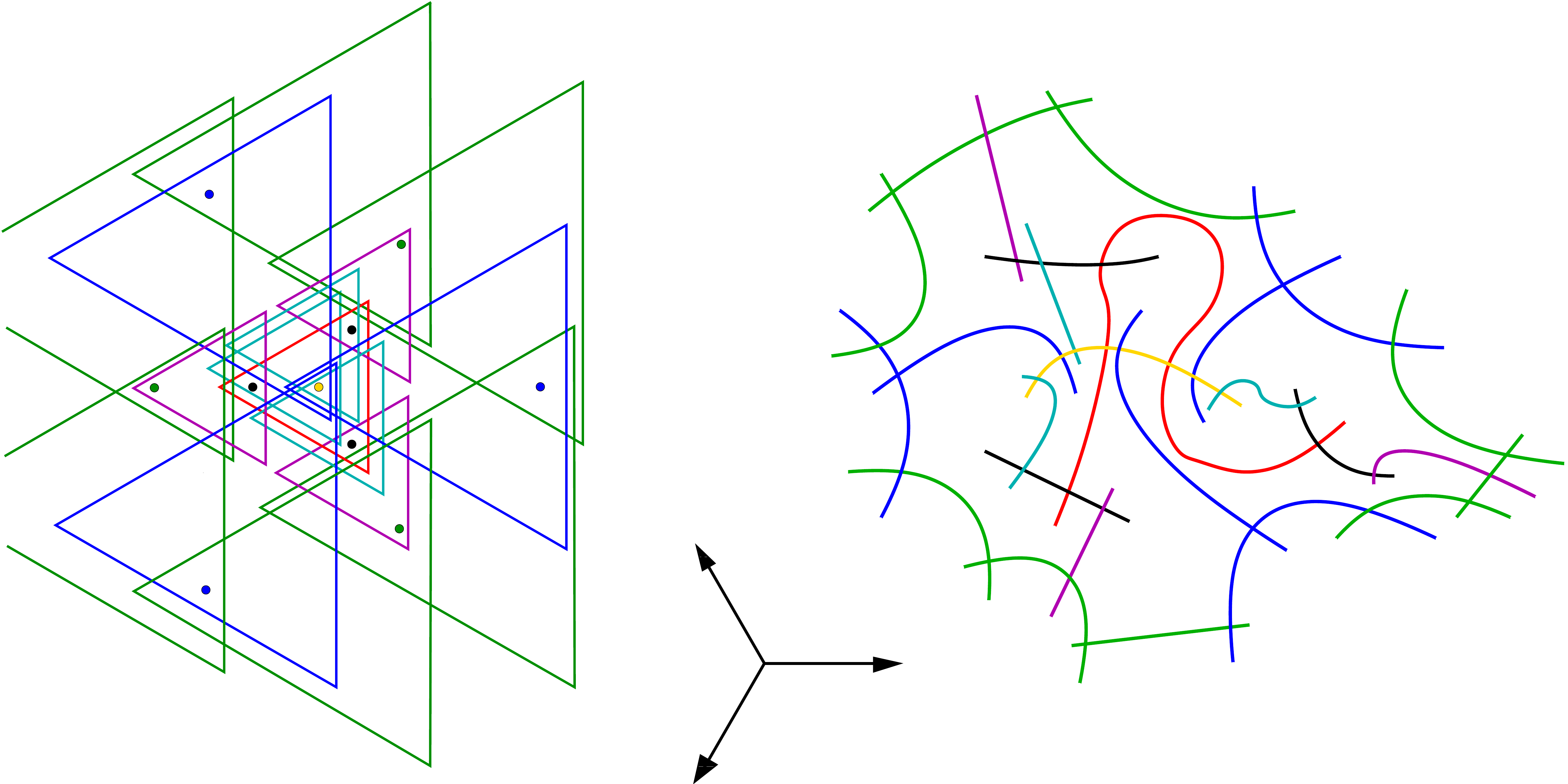}
    \caption{The drawing on the left defines an inclusion order of homothetic triangles. This height-$2$ order does not have a pseudosegment representation.\label{fig:drei-dim-nonPseudo-4}}
    \end{figure}%

%%%%%%%%%%%%%%%%%%%%%%%%%%%%%%%%%%%%%%%%%%%%%%%%%%%%%%%%

\Proof 
The left drawing in Figure~\ref{fig:drei-dim-nonPseudo-4} defines a poset $P$ by ordering the homothetic triangles by inclusion.
Some of the triangles are so small that we refer to them as points from now on.
Each inclusion in $P$ is witnessed by a point and a triangle, and hence $P$ has height~$2$.
To see that it is $3$-dimensional we use the drawing and the three directions depicted in Figure~\ref{fig:drei-dim-nonPseudo-4}.
By applying the same method as we did for Proposition~\ref{prop:GIGdim4} we obtain three linear extensions forming a realizer of $P$.

We claim that $P$ is not a pseudosegment intersection graph\footnote{The intersection graph of curves where each pair of curves intersects in at most one point.}, and hence not a GIG. 
Suppose to the contrary that it has a pseudosegment representation.
The six green triangles together with the three green and the three blue points form a cycle
of length $12$ in $G$.  Hence, the union of the corresponding pseudosegments in the representation contains a closed curve in
$\mathbb{R}^2$. Without loss of generality assume that the pseudosegment
representing the yellow point lies inside this closed curve (we may change the
outer face using a stereographic projection). The pseudosegments of the
three large blue triangles intersect the yellow pseudosegment and one blue
pseudosegment (corresponding to a blue point) each. The yellow and the blue
pseudosegments divide the interior of the closed curve into three regions. We
show that each of these regions contains one of the pseudosegments representing
black points. 

Each purple pseudosegment intersects the cycle in a point that is incident to one of the three bounded regions.
Now, each black pseudosegment intersects a purple one.
If such an intersection lies in the unbounded region, then the whole black pseudosegment is contained in this region.
This is not possible as for each of the black pseudosegments there is a blue pseudosegment representing a small blue triangle that connects it to the enclosed yellow pseudosegment without intersecting the cycle.
Thus, the three intersections of purple and black pseudosegments have to occur in the bounded regions, and in each of them one.
It follows that each of the three bounded regions contains one black pseudosegment.

Now, the red pseudosegment intersects each of the three black pseudosegments.
Since they lie in three different regions whose boundary it may only traverse
through the yellow pseudosegment, it has to intersect the yellow pseudosegment
twice. This contradicts the existence of a pseudosegment representation.  \qed

In the following we give constructions to show that

\begin{align*}
\textup{Stick}&\not\subset \textup{UGIG} &
\textup{UGIG}&\not\subset \textup{Stick}\\
\textup{BipHook}&\not\subset\textup{3-DORG} &
\textup{BipHook}&\not\subset \textup{Stick}
\end{align*}

\begin{figure}[htb]
    \centerline{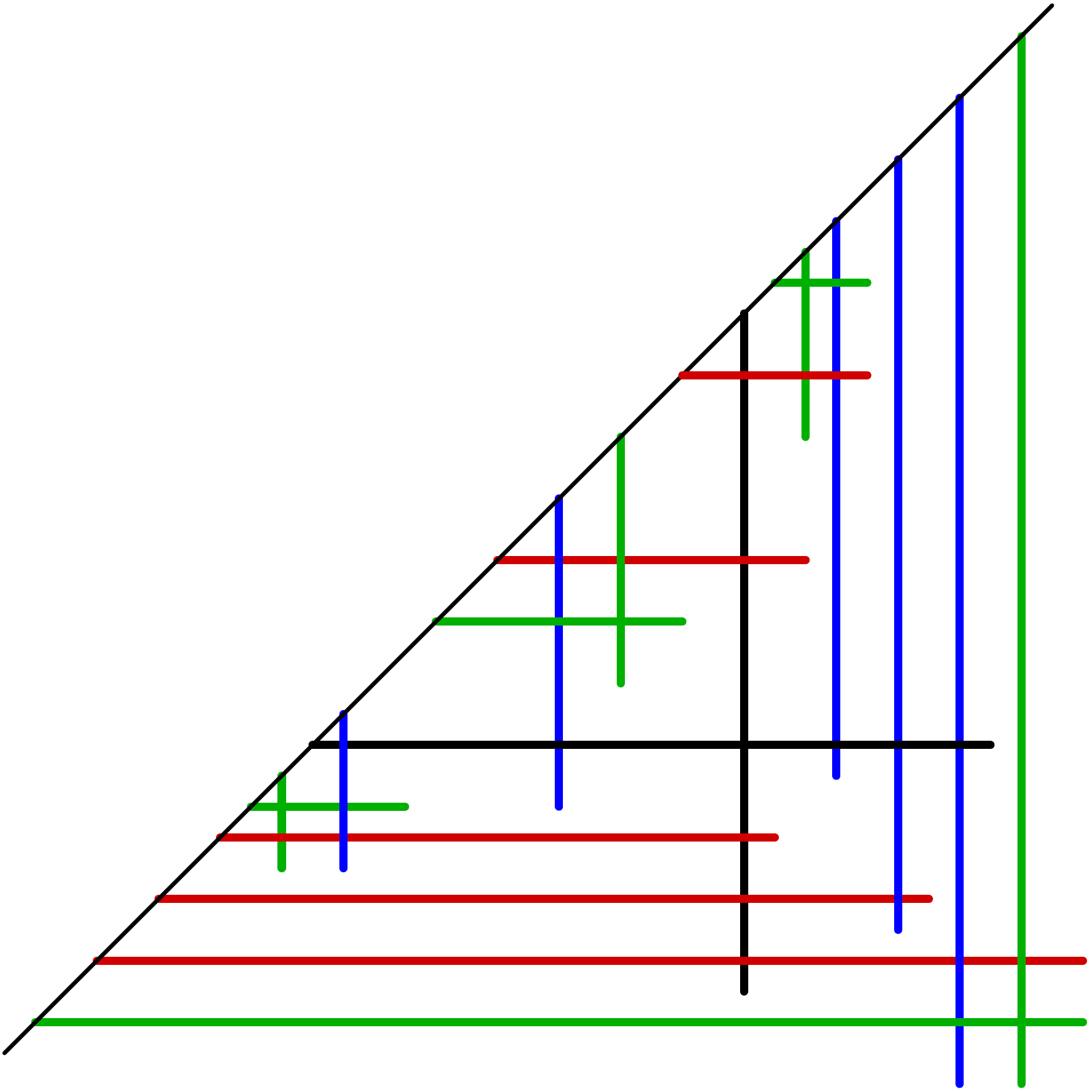}
    \caption{A stick representation of a graph that is not a
UGIG.\label{fig:STICKnonUGIG}}
    \end{figure}

\begin{proposition}\label{prop:STICKnotUGIG}
The Stick graph shown in Figure~\ref{fig:STICKnonUGIG} is not a UGIG.
\end{proposition}

%%%%%%%%%%%%%%%%%%%%%%%%%%
% in einem figure environment mit caption
\calc_figscale{30}%

%%%%%%%%%%%%%%%%%%%%%%%%%%
\Proof
Let $G$ be the graph represented in Figure~\ref{fig:STICKnonUGIG}.
Let $v$ and $h$ be the two adjacent vertices of $G$ that are drawn as black sticks in the figure.
There are five pairs of intersecting blue vertical and red horizontal segments
$v_1,h_1,\dots,v_5,h_5$.
Each $v_i$ intersects $h$ and each $h_i$ intersects $v$.
Four of the pairs $v_i,h_i$ form a 4-cycle with a pair of
green segments $q_i,r_i$.

\begin{figure}[htb]
    \centerline{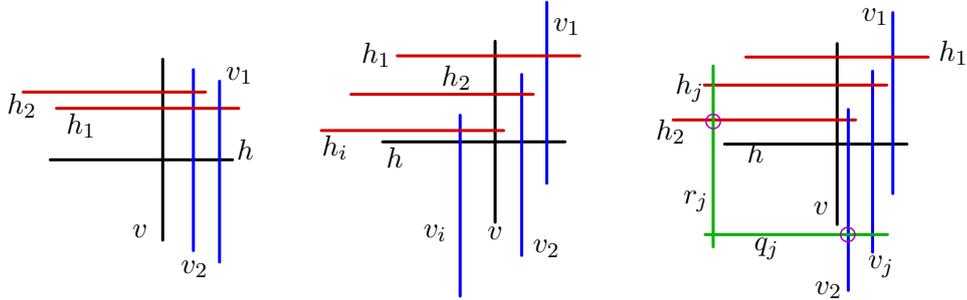}
    \caption{Left: The intersection points $p_1$, $p_2$ in the upper right quadrant form a chain in $<_\dom$. Middle: $p_i$ does not dominate $p_2$ in $<_\dom$.
Right: The green segments $q_j,r_j$ for the middle pair of segments $h_j,v_j$ cannot be
added.\label{fig:STICKnotUGIG}}
    \end{figure}

Suppose that $G$ has a UGIG representation.
We claim that in any such representation the intersection points $p_i$ of $v_i$ and $h_i$ form a chain in $<_\dom$ after a suitable rotation of the
representation.
Note that one quadrant formed by the segments $v$ and $h$ (without
loss of generality the upper right one) contains at least two
of the $p_i$'s by the pigeonhole principle. Assume without loss of generality
that $p_1$ and $p_2$ lie in this quadrant.
If $p_1$ and $p_2$ are incomparable in $<_\dom$, then the horizontal
segment $h_1$ of the lower intersection point has a forbidden intersection with the vertical segment $v_2$ of the higher one, see Figure~\ref{fig:STICKnotUGIG} left.

%%%%%%%%%%%%%%%%%%%%%%%%%
% in einem figure environment mit caption
\calc_figscale{30}%

%%%%%%%%%%%%%%%%%%%%%%%%%%

So $p_1$ and $p_2$ are comparable in $<_\dom$.
We may assume that $v_2,h_2$ is the pair of segments whose
intersection point is dominated in $<_\dom$ by all other intersection points in
the upper right quadrant. We observe that the lower endpoint of $v_2$ lies
below the lower endpoint of $v$, and the left endpoint of $h_2$ lies to the
left of the left endpoint of $h$ as shown in the middle of Figure~\ref{fig:STICKnotUGIG}.
It follows that if an intersection point $p_i$ does not dominate $p_2$, then $p_i$ lies below $h_2$ and to the left of $v_2$, but not in the upper right quadrant by our choice of $p_2$ (see Figure~\ref{fig:STICKnotUGIG} for an example).
It is easy to see that the remaining two intersection points $p_j$ ($j\not\in\{1,2,i\}$) then have to dominate $p_2$ in $<_\dom$, as otherwise we would see forbidden intersections among the blue and red segments.

We conclude that, in each case, four of the points $p_1,\dots, p_5$ lie in the upper
right quadrant and that they form a chain with respect to $<_\dom$.
Thus at least one pair of segments $v_j,h_j$ with $p_j$ being in the middle of the chain has neighbours $q_j,r_j$.
However, as indicated in the right of Figure~\ref{fig:STICKnotUGIG}, $q_j$ and $r_j$
cannot be added without introducing forbidden intersections.
Hence $G$ does not have a UGIG representation.
\qed

We now show that there is a 3-DORG that is not a BipHook graph.
We will use the following lemma for the argument.
\begin{lemma}\label{lem:twin_hook}
Let $G$ be a bipartite graph and $G'$ be the graph obtained by adding a twin to each vertex of $G$ (i.e., a vertex with the same neighbourhood). Then $G'$ is a hook graph if and only if $G$ is a Stick graph.
\end{lemma}
\Proof
Suppose that $G'$ has a hook representation.
Consider twins $v,v'\in V(G)$ and the position of their neighbours in a hook representation.
Suppose that there are vertices $u,w\in N(v)$, such that the order on the diagonal is $u,v,v',w$.
One can see that this order of centres together with edges $uw$ and $v'u$ would
force the hooks of $v$ and $v'$ to intersect, which contradicts their
non-adjacency.
Thus either no neighbour of $v$ occurs before $v$ or no neighbour of $v'$ 
occurs after $v'$ on the diagonal.
This shows that the hook of $v$ or $v'$ can be drawn as a stick, and it follows that $G$ has a Stick representation.

Conversely, in a stick representation of $G$ twins can easily be added to obtain a stick representation of $G'$.
\qed

\begin{figure}[htb]
    \centerline{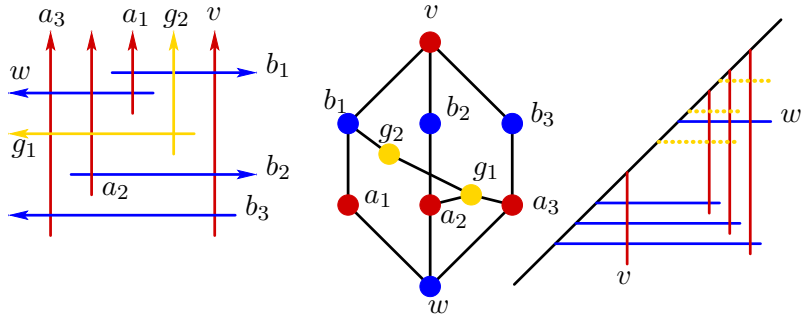}
    \caption{A \textup{3-DORG} that is not a Stick graph.\label{fig:3DORGnotHOOK}}
    \end{figure}

\begin{proposition}\label{prop:3dorgnothook}
The \textup{3-DORG} in Figure~\ref{fig:3DORGnotHOOK} is not a Stick graph.
\end{proposition}
%
%%%%%%%%%%%%%%%%%%%%%%%%%%%%%%%%%%%%%%%%%%%%%%%%%%%%%%%%%%%%
% in einem figure environment mit caption
\calc_figscale{27}%

%%%%%%%%%%%%%%%%%%%%%%%%%%%%%%%%%%%%%%%%%%%%%%%%%%%%%%%%%%%%
%
\Proof Suppose to the contrary that a Stick representation of the graph exists. We may assume
that $v$ is a vertical and $w$ a horizontal stick.  Observe that $w$ has to
lie above $v$ on the diagonal: Otherwise, two of the $a_i$'s have to lie either
before $v$ or after $v$, however, for the outer one of such a pair of $a_i$'s
it is impossible to place a stick for $b_i$ that also intersects $v$.  Hence,
the Stick representation of $v,w$ and the $a_i$'s and $b_i$'s have to look as
in Figure~\ref{fig:3DORGnotHOOK}.  By checking all possible positions of
$g_1$, i.e., permutations of $\{a_1,a_2,a_3\}$ and the correspondingly forced
permutation of $\{b_1,b_2,b_3\}$ in the representation, it can be verified that the
representation cannot be extended to a representation of the whole graph. The
cases are indicated in Figure~\ref{fig:3DORGnotHOOK}.
\qed

As a consequence, there is a \textup{3-DORG} that is not a bipartite hook graph.
Indeed, if we add a twin to each vertex of the graph shown in Figure~\ref{fig:3DORGnotHOOK} then the obtained graph is still a 3-DORG.
It can not be a BipHook graph as otherwise by Lemma~\ref{lem:twin_hook} we would conclude that the graph in Figure~\ref{fig:3DORGnotHOOK} is a Stick graph.

We next show a construction of a bipartite hook graph that is not a Stick graph.
A related construction was also presented in~\cite{tomMaster}.

%%%%%%%%%%%%%%%%%%%%%%%%%%%%%%%%%%%%%%%%%%%%%%%%%%%%%%%%%%%%
% in einem figure environment mit caption
\calc_figscale{55}%
\begin{figure}[htb]
    \centerline{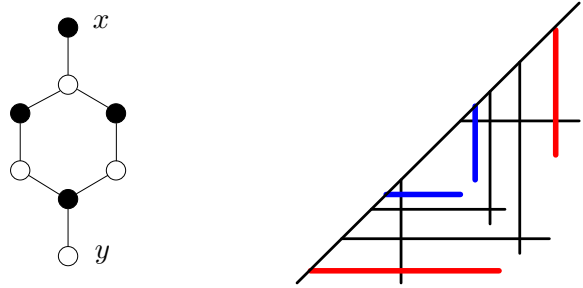}
    \caption{The graph $\Phi$ and the two possible positions of
  $x$ and $y$ in a Stick representation of $G$.\label{fig:nonstickhookG6}}
    \end{figure}%

%%%%%%%%%%%%%%%%%%%%%%%%%%%%%%%%%%%%%%%%%%%%%%%%%%%%%%%%%%%%

\begin{proposition}\label{prop:hooknotstick}
There is a bipartite hook graph that is not a Stick graph.
% $\mbox{BipHook}\not\subseteq\mbox{Stick}$.
\end{proposition}

\Proof The proof is based on the graph $\Phi$ shown in
Figure~\ref{fig:nonstickhookG6}. The vertices $x$ and $y$ are the {\em
  connectors} of $\Phi$. Let $G$ be a graph that contains
an induced $\Phi$ and a path $p_{xy}$ from $x$ to $y$ such that there is no adjacency between inner vertices of $p_{xy}$ and the $6$-cycle of $\Phi$. Observe that the Stick representation of
the 6-cycle is essentially unique. Now it is easy to check that in a Stick representation of $G$ the sticks for
the connectors have to be placed like the two blue sticks or like the two red
sticks in Figure~\ref{fig:nonstickhookG6}, otherwise the sticks of $x$ and $y$
would be separated by the $6$-cycle, whence one of the sticks representing
inner vertices of $p_{xy}$ and a stick of the 6-cycle would intersect. Depending on the placement the connectors are of type {\em inner} (blue) or
{\em outer} (red).

\begin{figure}[h]
 \centering
 \includegraphics{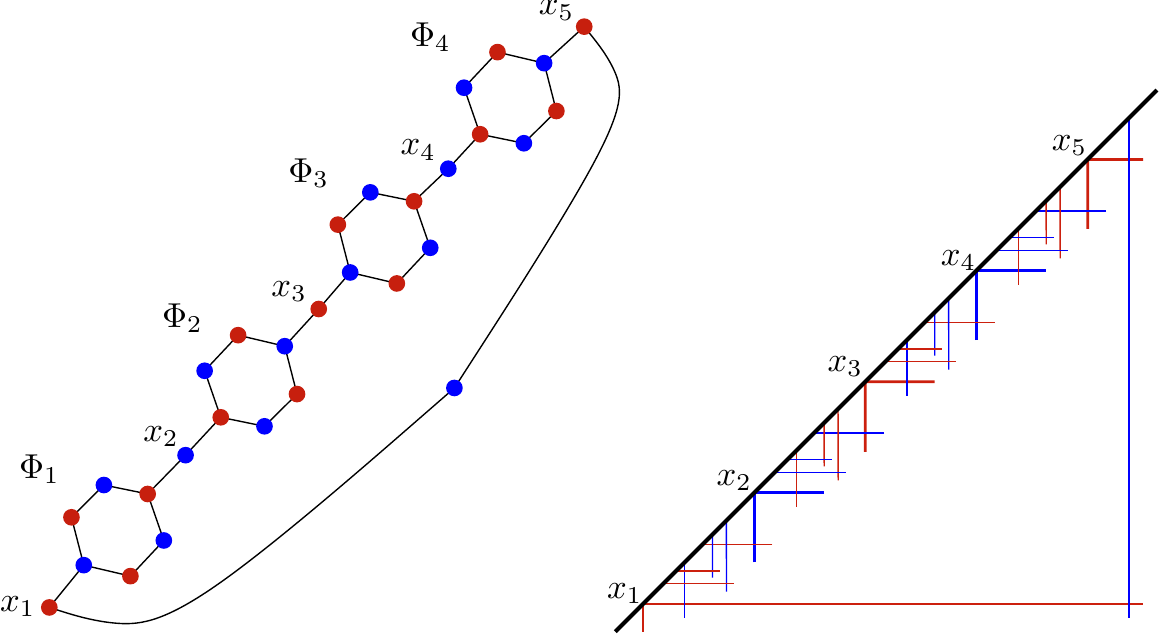}
 \caption{A bipartite hook graph (with hook representation) that is not a Stick graph.}
 \label{fig:nonstickhook}
\end{figure}

\begin{figure}[htb]
    \centerline{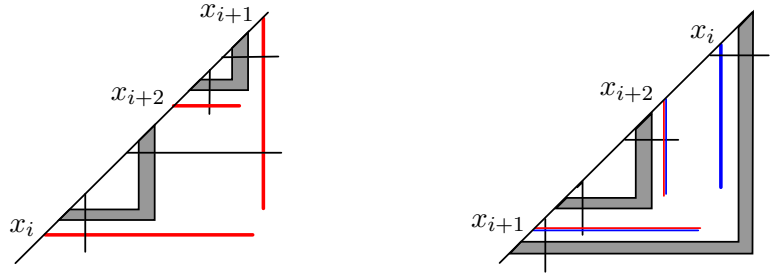}
    \caption{Stick representations of $\Phi_i$ and $\Phi_{i+1}$ with connecters of type inner--inner (left) and inner--outer (right).\label{fig:nonstickhookproof2}}
    \end{figure}

Consider the graph $\Phi^{4}$ depicted in Figure~\ref{fig:nonstickhook} together with a hook representation of it.
Suppose for contradiction that $\Phi^4$ has a Stick representation.
It contains four copies $\Phi_1,\ldots,\Phi_4$ of the graph $\Phi$ with connectors $x_1,\ldots,x_5$.
By our observation above, the connectors of each $\Phi_i$ are either of
inner or outer type.
We claim that for each $i\in\{1,2,3\}$, connecters of $\Phi_i$ and $\Phi_{i+1}$ are of different type.
If the type of both connectors of $\Phi_i$ and $\Phi_{i+1}$ is inner, then such a placement would force extra edges, specifically an edge between the two 6-cycles of $\Phi_i$ and $\Phi_{i+1}$.
And if both are outer then such a placement would separate $x_i$ and $x_{i+2}$, see
Figure~\ref{fig:nonstickhookproof2} on the left.

%%%%%%%%%%%%%%%%%%%%%%%%%%%%%%%%%%%%%%%%%%%%%%%%%%%%%%%%%%%%%%%%%%%%%%%%%%
% in einem figure environment mit caption
\calc_figscale{35}%

%%%%%%%%%%%%%%%%%%%%%%%%%%%%%%%%%%%%%%%%%%%%%%%%%%%%%%%%%%%%%%%%%%%%%%%%%%

It follows that the connector type of the $\Phi_i$'s is alternating.
In particular, there is $i\in\{1,2\}$ such that the connectors of $\Phi_{i},\Phi_{i+1},\Phi_{i+2}$ are of type inner--outer--inner in this order.
The right-hand side of Figure~\ref{fig:nonstickhookproof2} illustrates how $\Phi_{i}$ and $\Phi_{i+1}$ have to be drawn in a Stick representation.
Since $x_{i+2}$ is one of the inner type connectors of $\Phi_{i+2}$, there is no chance of adding the sticks for $\Phi_{i+2}$ to the drawing without intersecting sticks representing $\Phi_{i}$ and $\Phi_{i+1}$.
This is a contradiction and hence $\Phi^4$ is not a Stick graph.\qed

\subsection{Stabbability}\label{subsec:stab}
We proceed to show that
$$\textup{SegRay}\not\subset \textup{StabGIG}\qquad\qquad
\textup{4-DORG}\not\subset \textup{StabGIG}.$$

As an intermediate step we prove that there are GIGs that are not stabbable.
Techniques used in the proof will be helpful to show the two seperations.

\begin{proposition}\label{prop:nonStabGIG}
 There exists a GIG that is not a StabGIG.
\end{proposition}

\Proof
Consider a GIG representation of a complete bipartite graph $K_{n,n}$.
The GIG representation forms a grid in the plane. Now we add segments 
such that for every pair of cells in the same row or in the same column
there is a segment that has endpoints in both of the cells. Furthermore,
those segments can be drawn in such way that a horizontal and a vertical segment
intersect if and only if both intersect a common cell completely, that is, they do not have an endpoint in this cell.
Denote the resulting GIG representation by $R_n$ and the corresponding GIG by $G_n$. 

\begin{figure}[htb]
    \centerline{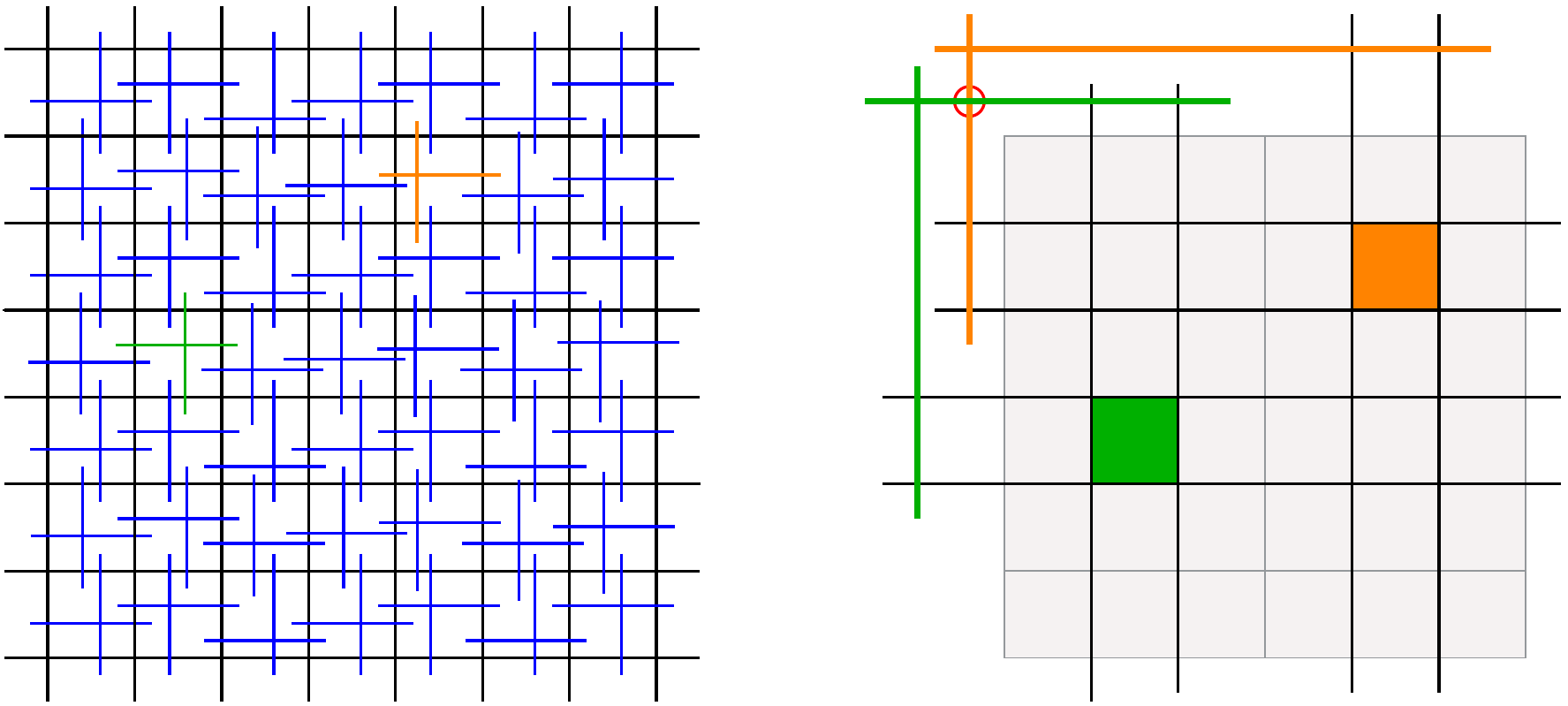}
    \caption{Partial representation of $R_n$ on the left. Replaced cell segments in $R'_n$ on the right.\label{fig:gignonstab}}
    \end{figure}

Suppose for contradiction that $G_n$ has a stabbable GIG representation $R'_n$ for all $n\in\mathbb{N}$.
By the Erd\H{o}s-Szekeres theorem for monotone subsequences, for every $k\in\mathbb{N}$ there exists $n\in \mathbb{N}$ such that in $R_n$ there are
subsets $H$ and $V$ consisting of $k$ horizontal and $k$ vertical segments that represent vertices of the $K_{n,n}$ in $G_n$, such that they appear in the same order (up to reflection) as segments in $R_n$ representing the same set of vertices.
In $R_n$ those segments induce a subgrid to which we added the blue segments depicted in Figure~\ref{fig:gignonstab}.
That is, for each cell $c$ in the subgrid we have a horizontal segment $h_c$ and a vertical one $v_c$ such that $h_c$ and $v_c$ intersect only each other and the segments building the boundary of $c$.

%%%%%%%%%%%%%%%%%%%%%%%%%%%%%%%%%%%%%%%%%%%%%%%%%%%%%%%%%%%%%%%%%%%%%%%%%%
% in einem figure environment mit caption
\calc_figscale{60}%

%%%%%%%%%%%%%%%%%%%%%%%%%%%%%%%%%%%%%%%%%%%%%%%%%%%%%%%%%%%%%%%%%%%%%%%%%%

It is easy to see that this partial representation in $R_n$ is not stabbable if $k$ is large.
Since the segments of the subgrid appear in the same order (up to reflection) in $R'_n$, we only have to consider the placement of cell segments $h_c$ and $v_c$.  
We restrict our attention to cells not lying on the boundary of the grid and fix a stabbing line $\ell$ for $R'_n$.
There are two possibilities for the placement of $h_c$ and $v_c$ in $R'_n$.
One case is that the intersection point $p_c$ of $h_c$ and
$v_c$ lies in $c$ or in one of the eight cells surrounding $c$. Then the segments
$h_c$ and $v_c$ can only be stabbed by $\ell$ if at least one of those eight cells around $c$ is
intersected by $\ell$.  The cells intersected by $\ell$ in $R'_n$ are only $O(k)$ many, so their
neighbouring cells are only $O(k)$ many as well.  This shows that $\Omega(k^2)$
intersection points $p_c$ have to lie outside of the grid in $R'_n$ (as depicted in Figure~\ref{fig:gignonstab} on the right). 
However, we show that this is possible for only $O(k)$ of them.

If an intersection point lies outside of the grid it is assigned to one quadrant,
i.e., $p_c$ lies above or below and left or right of the interior of the grid.
Every quadrant contains at most $O(k)$ points $p_c$: We index each cell by
its row and column in the grid so that the bottom- and leftmost cell is $c_{1,1}$.  If the intersection points corresponding to cells
$c_{u,v}$ and $c_{x,y}$ lie in the upper left quadrant, then $u<x$ implies
$y\leq v$.
This is illustrated in Figure~\ref{fig:gignonstab} where it is shown that otherwise the cell segments of the colored cells produce a forbidden intersection.
It follows that at most $O(k)$ intersection points of cell segments can lie in one quadrant, and hence $O(k)$ of them lie outside of the grid.
We conclude that $G_n$ has no
stabbable GIG representation for a sufficiently large $n$.
\qed

For SegRay graphs we give a similar construction that shows that there are
SegRay graphs which do not belong to StabGIG. First we will construct a graph
that cannot be stabbed in any SegRay representation.  
\begin{lemma}\label{lem:segrayunique}
  Let $R$ be a SegRay representation of a cycle $C$ with $2n$ vertices. 
  For the vertices in $C$ being represented as rays it holds that their order in $C$ is up to reflection and cyclic permutation equal to the order of the rays representing them in $R$ .
\end{lemma}
\Proof Let $L$ be a horizontal line below all horizontal segments in $R$.
Contracting each ray to its intersection point with $L$ yields a planar drawing of a cycle $C'$ with $n$ vertices such that the vertices lie on $L$ and edges are drawn above $L$.
This is also known as a \emph{1-page embedding} of $C'$.
It is easy to see that edges in $C'$ have to connect consecutive vertices on $L$ or the two extremal ones.
Now the conclusion of the lemma is straightforward.
\qed

\begin{proposition}\label{prop:nonstab-segray-rep}
There exists a SegRay graph that has no stabbable SegRay representation.
\end{proposition}

%%%%%%%%%%%%%%%%%%%%%%%%%%%%%%%%%%%%%%%%%%%%%%%%%%%%%%%%%%
% in einem figure environment mit caption
\calc_figscale{30}%
\begin{figure}[htb]
    \centerline{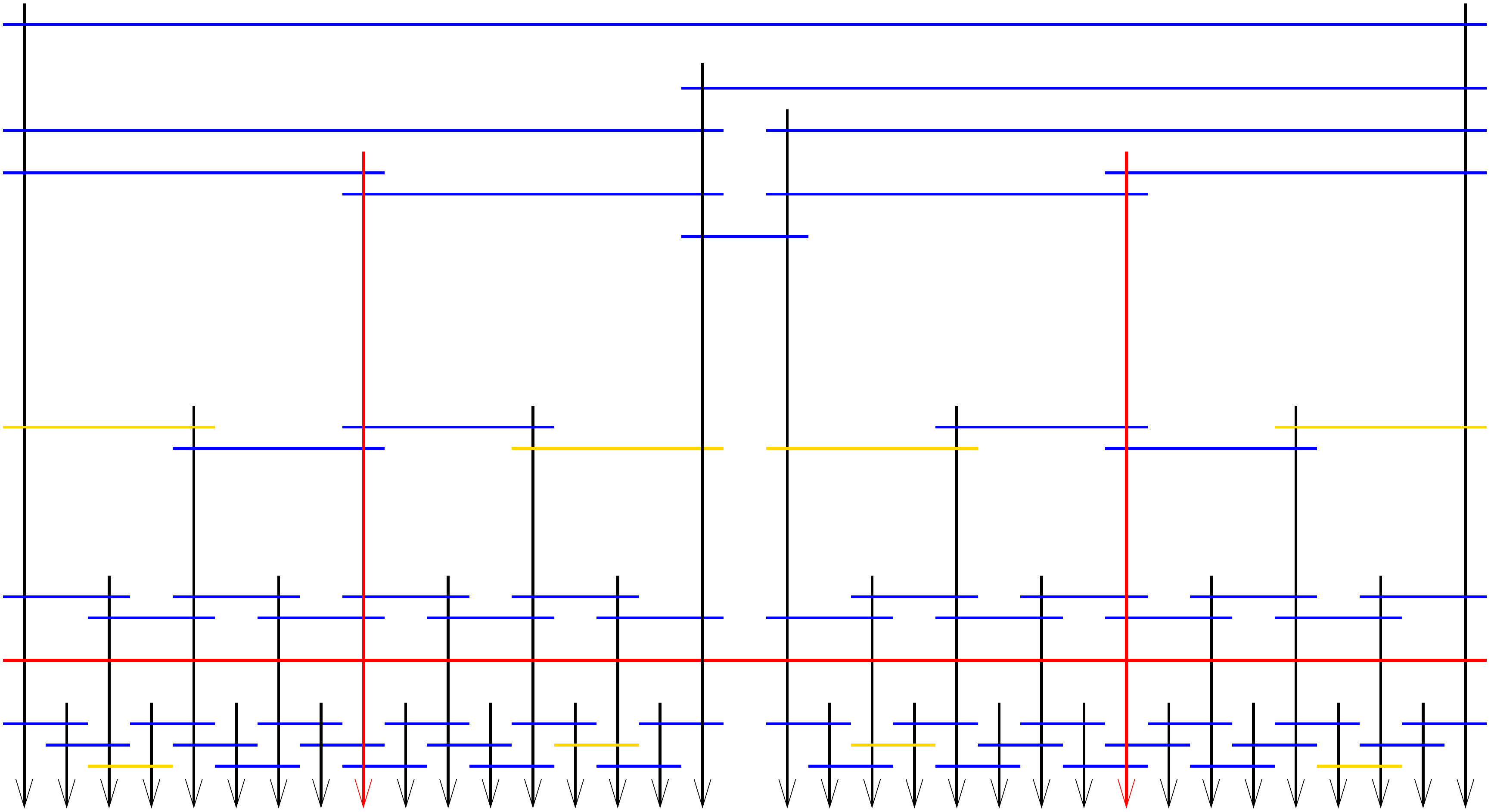}
    \caption{A SegRay graph with no stabbable SegRay representation.\label{fig:nonstabsegray}}
    \end{figure}%

%%%%%%%%%%%%%%%%%%%%%%%%%%%%%%%%%%%%%%%%%%%%%%%%%%%%%%%%%%

\Proof
Consider the graph defined by the SegRay representation $R$ in Figure~\ref{fig:nonstabsegray}.
Let $R'$ be an arbitrary SegRay representation of this graph.
The order of the rays in $R'$ is up to reflection and cyclic permutation equal to the one in $R$ by Lemma~\ref{lem:segrayunique}.
Hence without loss of generality the rays of the left half in $R$ appear consecutively in $R'$.
Now observe that the two yellow segments below the red segment in the left half of $R$ also have to be below the red segment in $R'$.
Similarly, the two yellow segments above the red segment in $R$ must lie
above the red segment in $R'$.  
Furthermore, these two yellow segments have to lie below the top of the red vertical ray in $R'$.
It follows that, as in $R$, the red segment and the red ray seperate the plane into four quadrants in $R'$ such that each quadrant contains exactly one of the considered yellow segments.
Any line in the plane can intersect at most three of the quadrants and thus will miss a yellow segment in $R'$.
Therefore, $R'$ is not stabbable and the conlusion follows.  
\qed

We add a vertex $h$ to the graph in Figure~\ref{fig:nonstabsegray} that
is adjacent to all rays.
This graph is still a SegRay graph.
We call this graph  a \emph{bundle} and the set of horizontal segments its \emph{head}.
The bundle is not stabbable in any SegRay representation by the proposition above.
This means, in any StabGIG representation of a bundle there is one horizontal
segment above the segment representing $h$ and one below.
Indeed, otherwise the representation of the bundle can be modified by extending vertical segments to rays in one direction to obtain a stabbed SegRay representation.
Using this property of a bundle we can show the following.

\begin{figure}[b]
 \centering
 \includegraphics{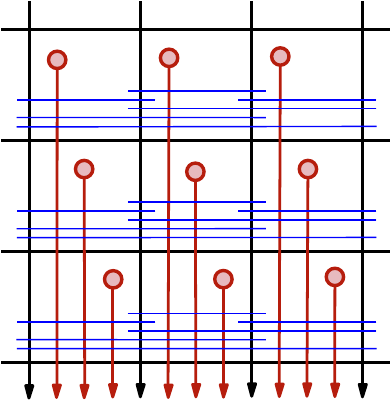}
 \caption{Illustration of a SegRay graph that is not a StabGIG}
 \label{fig:nonstabgig_segray}
\end{figure}

\begin{proposition}\label{prop:nonstabsegray}
  There exists a SegRay graph that is not a stabbable GIG.
\end{proposition}

\Proof Similar to the construction in Proposition~\ref{prop:nonStabGIG},
consider a SegRay representation of a complete bipartite graph $K_{n,n}$.
In this representation we see a grid with cells.
We place in each of the cells the head of a bundle as indicated in Figure~\ref{fig:nonstabgig_segray}. Now,
for each pair of cells in the same row of the grid, add a spanning horizontal
segment
with endpoints in the given cells.
We do it in such a way that the rays of a bundle are intersected by the segment if the head of the bundle lies in a cell between the two given cells.  

Denote by $R_n$ the resulting SegRay representation and let $H_n$ be the SegRay graph defined by $R_n$.
Suppose that $H_n$ has a stabbable GIG representation $R'_n$.
As in the proof of Proposition~\ref{prop:nonStabGIG}, given an integer $k\geq 1$ it follows by the Erd\H{o}s-Szekeres theorem that for sufficiently large $n$ there is a subgrid of size $k$ in $R'_n$, where the order of the horizontal and vertical segments is either preserved or reflected with respect to $R_n$.
Assume that it is preserved.
Now we restrict our attention to the relevant bundles and horizontal segments of $R_n$ according to the cells of the subgrid.
In $R_n$ again this looks like in Figure~\ref{fig:nonstabgig_segray}, but this time with respect to the fixed subgrid.

\begin{figure}[h]
 \centering
 \includegraphics{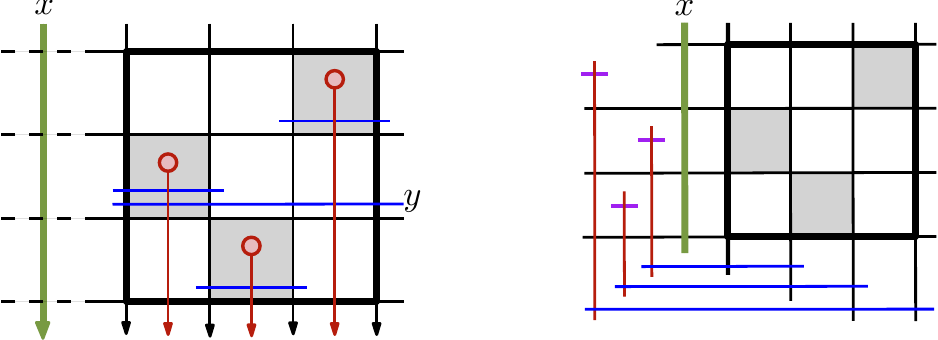}
 \caption{For at least one grey cell $c$, the purple segment $h_c$ lies to the right of $x$ in $R'_n$.}
 \label{fig:3x3-box}
\end{figure}

Let us now consider the placement of the bundles and blue segments in $R'_n$.
Given a cell $c$ in the subgrid, let $y_c$ be the horizontal grid segment bounding $c$ from below.
By Proposition~\ref{prop:nonstab-segray-rep} and its consequences, the bundle lying in cell $c$ contains a horizontal segment that lies above $y_c$ in $R'_n$.
We denote this segment by $h_c$ and let $x_c$ be an arbitrary ray of the bundle intersecting $h_c$.
Consider now the left side of Figure~\ref{fig:3x3-box} showing a $3\times 3$ box and a ray $x$ of the subgrid that lies strictly to the left of the box in the representation $R_n$.
Let $c_1,c_2,c_3$ be the three shaded cells.
Then we claim that at least one of $h_{c_1},h_{c_2},h_{c_3}$ is placed to the right of $x$ in $R'_n$.

Suppose to the contrary that all lie to the left.
If we use the fact that $h_{c_i}$ is above $x_{c_i}$ in $R'_n$ for each $i\in \{1,2,3\}$, then it is straightforward to see that $h_{c_1},h_{c_2},h_{c_3},x_{c_1},x_{c_2},x_{c_3}$ and the three short blue horizontal segments depicted on the left of Figure~\ref{fig:3x3-box} have to be placed in $R'_n$ as shown on the right of Figure~\ref{fig:3x3-box} (segments $h_{c_i}$ are colored purple).
But then the segment $y$, which is the long blue one on the left, cannot be added to the partial representation without creating forbidden crossings.
This shows our subclaim.

\begin{figure}[h]
 \centering
 \includegraphics{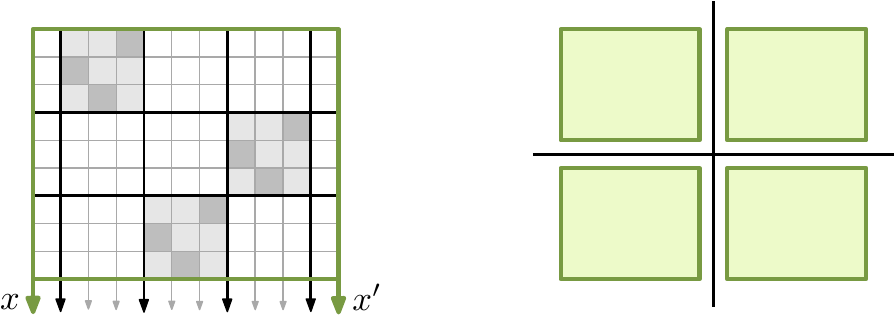}
 \caption{There is a cell $c$ in the green box such that $h_c$ is drawn inside the box in $R'_n$.}
 \label{fig:9x9-box}
\end{figure}

In the next step we consider the green box of the fixed subgrid shown on the left of Figure~\ref{fig:9x9-box}.
Using our subclaim we have that each of the three shaded $3\times 3$ boxes contains a cell $c$ such that $h_c$ is placed to the right of $x$ in $R'_n$.
Now we apply the symmetric version of this claim to deduce that one of these three segments also lies to the left of $x'$ in $R'_n$.
We conclude that there is a segment that is strictly contained in the green box in $R'_n$.

In the final step we consider four copies of the green box that are placed in the fixed subgrid of $R_n$ as shown on the right of Figure~\ref{fig:9x9-box}.
Since each copy strictly contains a segment in $R'_n$, each line in the plane will miss at least one of the four segments.
This shows that $R'_n$ is not stabbable for sufficiently large $n$ and completes the proof.
\qed

\begin{proposition}\label{prop:4dorgnotstabbable}
There exists a 4-DORG that is not a stabbable in any GIG representation.
\end{proposition}

\Proof
 Since the ideas here are similar to those used for Propositions~\ref{prop:nonstab-segray-rep} and \ref{prop:nonstabsegray}, we only provide a sketch of this proof.
 Consider the following construction.
 Take a $4$-DORG representation of a complete bipartite
 graph $K_{n,n}$.
 Similarly to previous constructions this yields a grid with cells. 
 For each cell we add four rays starting in this cell, one in each direction and such that a vertical and a horizontal ray intersect if and only if they entirely intersect a common cell.
 Call this representation $R_n$ and the corresponding intersection graph $G_n$.
 We claim that for sufficiently large $n$ there is no stabbable GIG representation of $G_n$.

 Suppose to the contrary that there exists a StabGIG representation $R'_n$ of $G_n$.
 Again by applying the Erd\H{o}s-Szekeres theorem we may assume that there is a large subgrid of size $k$ in $R'_n$, such that the order of the grid segments in $R'_n$ agrees with the order in $R_n$ (up to reflection).
 
\begin{figure}[h]
 \centering
 \includegraphics{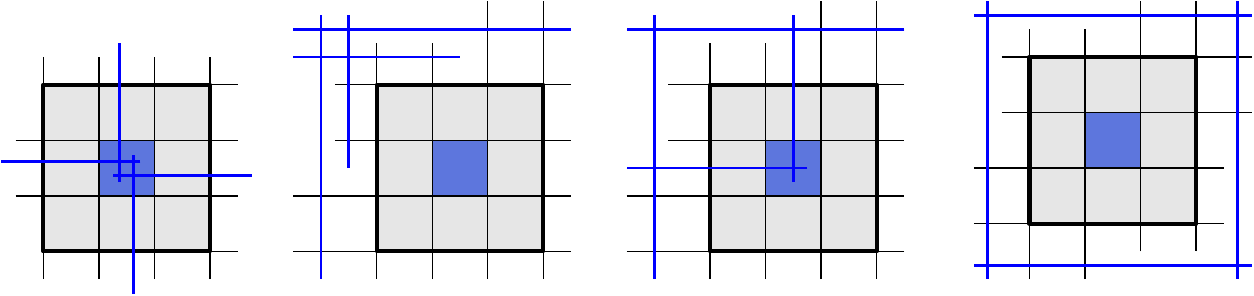}
 \caption{Four different situations for the blue cell and its segments in $R'_n$}
 \label{fig:nonstabdorg}
\end{figure}

 Given the representation $R'_n$, we want to partition the cells of the subgrid according to the placement of the four segments representing the rays that start in a given cell of our construction.
 Note that these four segments intersect in such a way that they enclose a rectangle in $R'_n$.
 Therefore, we can distinguish the following cases: the rectangle (1) is contained in a grid cell, (2) it does not intersect a grid cell, (3) it contains some but not all grid cells, and (4) it contains all of the grid cells (see Figure~\ref{fig:nonstabdorg} for the cases from left to right).
 Each of these cases again can be split into at most four natural subcases.
 For instance, if the rectangle contains some of the grid cells, then it also has to contain a corner of the grid, which gives rise to four subcases.
 
 Using similar arguments as in previous proofs of the paper and the assumption that $R'_n$ is stabbed by a line, one can show that each partion class contains at most $O(k)$ cells.
 Thus, for large enough $n$ and $k$ we get a contradiction since our subgrid has $k^2$ cells.
 This observation completes the proof.
\qed 

It remains open whether there exists a $3$-dimensional GIG that is not stabbable.
It is tempting to look for an example that produces again a large grid in every representation (to get non-stabbability), but it turned out that all these examples seem to have dimension $4$.
We also tried with SegRay graphs satisfying the properties of Lemma~\ref{lem:segrayprop} since they have dimension $3$.
However, we didn't succeed with finding such a SegRay that is not stabbable.

\section{Conclusion}

We have shown that Figure~\ref{fig:classinclusion} provides the correct inclusion
order of the given subclasses of GIG. An overview of the separating examples is
given in Table~\ref{tab:separating}.

The notion of order dimension was helpful in particular to exhibit examples
that separate classes. As a byproduct we have new insights regarding the interval
dimension of vertex-face posets of outerplanar maps (Corollary~\ref{cor:idim-opl}).

Another direction of research regarding these graph classes is recognition.
Currently the recognition complexity of some of the graph classes remains
open, see the table below.  We hope that our results help bringing these open
problems closer to a solution.

\begin{center}
\begin{tabular}{p{4cm}p{4cm}p{2cm}}
Class                        & recognition complexity  & reference \\[2pt] \hline\hline
GIG                          & NP-complete             & \cite{kratochvil1994} \\ 
UGIG                         & NP-complete             & \cite{mustata2013unit}\\
3-dim BipG                    & NP-complete             & \cite{felsner2015dim3h2} \\ \hline

3-dim GIG                    & Open                    &         \\
StabGIG                      & Open                    & \\
SegRay                       & Open                    &           \\
BipHook                      & Open                    &  \\
Stick                        & Open                    & \\ 
4-DORG                       & Open                    & \\
3-DORG                       & Open                    & \\ \hline

2-DORG                       & Polynomial              &\cite{STU10},\cite{cogis1982ferrers}\\ 
bipartite permutation        & Polynomial              & \cite{dm-pos-41} \\ 
\end{tabular}
\end{center}

%%%%%%%%%%%%%%%%%%%%%%%%%%%%%%%%%%%%%%%%%%%%%%%%%%%%%%%%%%%
\small
\advance\bibitemsep-1pt
\bibliographystyle{my-siam}
\bibliography{references}
%%%%%%%%%%%%%%%%%%%%%%%%%%%%%%%%%%%%%%%%%%%%%%%%%%%%%%%%%%%

\end{document}